\documentclass{article}
\usepackage[a4paper,left=1in, right=1in, top=1in, bottom=1in]{geometry}

\usepackage{multirow}
\newcommand{\set}[1]{\mathcal{#1}}
\usepackage{hyperref}
\usepackage{subfig}
\usepackage{placeins}
\usepackage{amsmath}
\usepackage{amsfonts}
\usepackage{amssymb}
\usepackage{graphicx}
\usepackage{adjustbox}
\usepackage{algorithm}
\usepackage[noEnd=false]{algpseudocodex}

\makeatletter
\newcommand{\algmargin}{\the\ALG@thistlm}
\makeatother
\algnewcommand{\StateLong}[1]{\State
    \parbox[t]{\dimexpr\linewidth-\algmargin}{\strut\hangindent=1em \hangafter=1 #1\strut}}

\usepackage{natbib}
 \bibpunct[, ]{(}{)}{,}{a}{}{,}

\date{}
\title{The Stochastic Dynamic Post-Disaster Inventory Allocation Problem with Trucks and UAVs}
\author{Robert van Steenbergen \cr Wouter van Heeswijk \cr Martijn Mes \cr \textit{University of Twente}}
\begin{document}
\maketitle

\begin{abstract}
Humanitarian logistics operations face increasing difficulties due to rising demands for aid in disaster areas. This paper investigates the dynamic allocation of scarce relief supplies across multiple affected districts over time. It introduces a novel stochastic dynamic post-disaster inventory allocation problem (SDPDIAP) with trucks and unmanned aerial vehicles (UAVs) delivering relief goods under uncertain supply and demand. The relevance of this humanitarian logistics problem lies in the importance of considering the inter-temporal social impact of deliveries. We achieve this by incorporating deprivation costs when allocating scarce supplies. Furthermore, we consider the inherent uncertainties of disaster areas and the potential use of cargo UAVs to enhance operational efficiency. This study proposes two anticipatory solution methods based on approximate dynamic programming, specifically decomposed linear value function approximation (DL-VFA) and neural network value function approximation (NN-VFA) to effectively manage uncertainties in the dynamic allocation process. We compare DL-VFA and NN-VFA with various state-of-the-art methods (e.g., exact re-optimization, PPO) and results show a 6-8\% improvement compared to the best benchmarks. NN-VFA provides the best performance and captures nonlinearities in the problem, whereas DL-VFA shows excellent scalability against a minor performance loss. From a practical standpoint, the experiments reveal that consideration of deprivation costs results in improved allocation of scarce supplies both across affected districts and over time. Finally, results show that deploying UAVs can play a crucial role in the allocation of relief goods, especially in the first stages after a disaster. The use of UAVs reduces transportation- and deprivation costs together by 16-20\% and reduces maximum deprivation times by 19-40\%, while maintaining similar levels of demand coverage, showcasing efficient and effective operations.
\end{abstract}

\section{Introduction}
\label{sec:introduction}
The Global Humanitarian Overview of the United Nations reports a record of 339 million people that needed humanitarian assistance and protection in 2023, which is a shocking increase compared to the 136 million people just five years earlier in 2018. Considering the increasing frequency and impact of disasters, humanitarian logistics is more important than ever. The purpose of humanitarian logistics is to provide affected people with supplies that are necessary for their survival, and when possible alleviate the human suffering caused by a lack of goods and services \citep{shao2020research}. Typically, humanitarian operations face a surge in demand directly in the response phase after a disaster, overwhelming the abilities of relief agencies \citep{zhou2017multi} and creating a significant scarcity of supplies \citep{najafi2013multi}. Given this scarcity, it is essential to utilize the available resources as effectively as possible \citep{holguin2012unique}. In the response phase, which generally ranges from a few days to a few weeks \citep{altay2006or}, humanitarian agencies face a spatio-temporal decision problem, requiring effective distribution of available supplies to multiple districts in a disaster area over time, in which they must consider the well-being of all affected people, the operational costs, and the limited funds that are received from donors \citep{huang2012models, huang2019equitable}.  

The duration of deprivation strongly impacts people's well-being. Most individuals can cope with short-term shortages, but long-term shortages are significantly more difficult to bear. As the timing of delivery matters -- rather than just the total amount delivered -- the inter-temporal impact of alleviating suffering due to delivery actions should be included in decision-making. For this purpose, \cite{holguin2013appropriate} defined the econometric concept of deprivation costs as the estimated value of human suffering caused by a lack of goods or services. The authors argue that the delivery effects must be assessed with multi-period models, as single-period formulations simply cannot account for inter-temporal effects. Due to the exponential growth of deprivation costs over time, decision makers are encouraged to focus on minimizing deprivation costs in situations with scarcity, whereas, in situations with sufficient supplies, the focus leans towards minimizing transportation costs \citep{perez2016inventory}. Incorporating deprivation costs in the objective could balance the conflicting goals and result in an operation that is both efficient in terms of low transportation costs and effective in terms of low deprivation costs.

While multi-period planning is important, anticipatory approaches are needed to tackle the inherent uncertainties associated with dynamically changing post-disaster areas \citep{hoyos2015or, besiou2020humanitarian}. Demand in response phases is uncertain and may evolve quickly due to environmental dynamics such as the mobilization of the population \citep{anaya2016models}. Agencies have to deal with inaccurate and evolving information regarding supply and demand. Obtaining correct estimates after disasters remains challenging \citep{lin2011logistics}, especially when reliable data for full planning horizons is unavailable \citep{sheu2010dynamic}. These uncertainties and dynamics prohibit the use of a fixed allocation plan for an entire response operation \citep{rottkemper2012transshipment, huang2015modeling}. 

Besides dynamic planning, new technologies also offer opportunities to deal with the evolving nature of a disaster response phase \citep{besiou2020humanitarian}. Humanitarian unmanned aerial vehicles (UAVs) exemplify such a new technology that could make humanitarian operations more effective and efficient. Cargo UAVs are flexible due to their speed, autonomous control, and ability to visit hard-to-reach places \citep{rejeb2021humanitarian, steenbergen2023}. Their flexibility and ability to deliver goods to many different locations offers the potential to address the spatio-temporal challenge of allocating scarce goods, while also adapting to the dynamically evolving disaster environment. Deploying multiple transportation modes has already yielded valuable results in humanitarian operations models, increasing cost-efficiency \citep{vanajakumari2016integrated}, decreasing response times \citep{alem2016stochastic}, improving security \citep{van2023heterogeneous}, and mitigating travel time uncertainty \citep{steenbergen2023}.

We established that (i) incorporating deprivation costs in a multi-period perspective is crucial to accurately measure the impact of humanitarian operations and to solve the intrinsic spatio-temporal decision problem, (ii) uncertainty and changing conditions are inherent to disaster areas and should be incorporated through anticipatory planning, and (iii) UAVs could play a key role in this setting. In this paper, we address these challenges and opportunities. We introduce the stochastic dynamic post-disaster inventory allocation problem (SDPDIAP) with uncertainty in supply and demand. The aim is to effectively allocate scarce supplies to multiple affected districts over time (i.e., spatio-temporal allocation) using trucks and UAVs, in a way that minimizes transportation costs and deprivation costs. The contributions of this work are (i) the formulation of a Markov decision process (MDP) model that captures the stochastic dynamic nature of the problem, (ii) the incorporation of deprivation costs in this stochastic dynamic problem, (iii) the development of a computationally efficient anticipatory solution method that handles the uncertainty and exploits the spatio-temporal nature of the problem by learning value functions for each individual district at each point in time, (iv) the development of a second solution method that captures the spatio-temporal and nonlinear problem elements with a neural network value function approximation, (v) numerical experiments that offer insights into the stochastic dynamic allocation of scarce supplies, trade-offs between transportation costs and deprivation costs, and the performance of our proposed methods, and (vi) the analysis of the opportunities and limitations of UAVs as an additional transportation mode for the dynamic allocation of relief supplies.

The remainder of this paper is structured as follows. Related literature on the stochastic dynamic post-disaster inventory allocation problem is reviewed in Section~\ref{sec:literature}. In Section~\ref{sec:model}, we describe the SDPDIAP and formulate the problem as an MDP. Section~\ref{sec:solution} introduces our solution methods and in Section~\ref{sec:experiments}, we describe the experimental setup and benchmarks to analyze the performance of our approaches. Section~\ref{sec:results} discusses the results and insights of the numerical experiments. Finally, the conclusions and future research directions are stated in Section~\ref{sec:conclusions}.

\section{Related Literature}
\label{sec:literature}
In this section, we cover related work regarding multi-period allocation problems in humanitarian logistics. We review works regarding stochastic dynamic supply allocation in humanitarian operations in Section \ref{sec:lit:sd-allocation}, reinforcement learning approaches in humanitarian resource allocation problems in Section \ref{sec:lit:rl-allocation}, and the variety of objective functions in humanitarian multi-period problems in Section \ref{sec:lit:objectives}. We conclude with the contribution statement in Section \ref{sec:lit:contributions}.

\subsection{Stochastic Dynamic Humanitarian Supply Allocation}
\label{sec:lit:sd-allocation}
To anticipate uncertainties and adapt to changes in a disaster area, various stochastic dynamic solution methods have been developed. Adjustable robust optimization is considered by \cite{ben2011robust} for allocating response operations and guiding evacuation traffic flows using a polyhedral uncertainty set for the demand. In a dynamic setting, the robust solution outperforms both deterministic- and stochastic programming counterparts. \cite{vanajakumari2016integrated} present an integrated location-allocation-routing model under budget constraints to minimize response times. New plans are dynamically generated to include updates of the realized demand and forecasted demand each period. Analysis of the results shows that a mix of different truck types is more cost-efficient than a single truck type. \cite{natarajan2017multi} study an inventory allocation problem for allocating scarce medical supplies to patients by formulating a finite horizon stochastic dynamic programming model. Their developed expected value heuristic and threshold heuristic obtain on average a 5\% gap against perfect information lower bounds and at least 10\% improvement compared to a first-come-first-serve heuristic that is often used in practice. \cite{liu2019robust} present a multi-period distribution model with helicopters that considers both injured people and relief commodities, minimizing the unsatisfied demand. Demand and supply are uncertain and a robust model predictive control approach is proposed to obtain robust plans and update them with new information. The study shows that the adjustments improve the results compared to a robust plan created for the entire horizon and the authors suggest an extension with multiple transportation modes for a more universal approach. \cite{zhan2021disaster} propose a dynamic multi-period relief allocation model with demand uncertainty. Sequential decisions are made using particle swarm optimization with a voting strategy. Considering limited supplies, decisions include serving a demand node now or temporarily neglecting the node. The objective is to minimize the total unmet demands weighted by a time-based penalty to balance between the egalitarian and utilitarian perspectives. \cite{ghasbeh2022equitable} develop a robust, multi-commodity, multi-stage model to address the facility location and relief distribution considering uncertainty in demand and budget. By sequentially solving the model using adaptive large neighborhood search (ALNS) and multi-dimensional local search (MDLS), robust distribution plans are obtained that minimize the unsatisfied demand and travel time. In general, the stochastic dynamic models are often not computationally efficient, hence various works have proposed (meta)heuristics to solve the problems each period with new information. A more generic approach for these sequential decisions is reinforcement learning, which is discussed in the section below.

\subsection{Reinforcement Learning for Humanitarian Allocation Problems}
\label{sec:lit:rl-allocation}
A variety of reinforcement learning methods, in which a policy is trained to make decisions, have been proposed to solve humanitarian allocation problems. \cite{yu2019rollout} present a rollout algorithm for a deterministic multi-period allocation problem in which for each period a single demand unit is allocated to one district. The aim is to minimize transportation- and deprivation costs. The algorithm provides competitive results against exact dynamic programming. \cite{yu2021reinforcement} further improved results with tabular Q-learning. Q-learning is compared with exact dynamic programming and a greedy heuristic. With deterministic instances up to 8 demand points and 12 periods, Q-learning yields optimal solutions with computation times up to six times smaller than exact dynamic programming. \cite{fan2022dhl} formulate an MDP for an allocation problem to minimize transportation costs and deprivation costs. With problem sizes up to 10 districts and 3 allocation units each period, deep Q-learning performs better than tabular Q-learning in terms of calculation time and the objective value, whereas an exact algorithm fails to provide an optimal solution within a reasonable time limit. A related application for reinforcement learning in humanitarian operations is the allocation of response teams to specific districts. \cite{nadi2016reinforcement} formulate an MDP model and apply tabular Q-learning to allocate teams to specific arcs such that the total assessment time is minimized. \cite{nadi2017adaptive} extend this work by formulating an MDP to allocate both assessment teams and response teams to different locations. The problem is solved with tabular Q-learning, considering uncertain demand and road damage. While various learning-based approaches have been proposed, they mostly focus on small-scale deterministic instances for which Q-values can be enumerated.

\subsection{Objectives in Humanitarian Operations}
\label{sec:lit:objectives}
Humanitarian operations often involve multiple objectives that need to be balanced. Besides transportation costs or travel distances, research focuses on response times \citep[e.g.,][]{najafi2014dynamic, lu2016real, vanajakumari2016integrated}, unsatisfied demands \citep[e.g.,][]{balcik2008last, lin2011logistics, liu2019robust, yang2023distributionally}, and equity \citep[e.g.,][]{tzeng2007multi, huang2015modeling, wang2023model}. \cite{vries2020optimization} argue that these kinds of objectives are more important than minimizing transportation costs when (i) resources are scarce and (ii) the urgency of incoming requests can be identified properly. In recent years, research focus has shifted towards including deprivation costs in the objective function \citep[e.g.,][]{rivera2016dynamic, perez2016inventory, moreno2018effective, huang2019equitable, yu2021reinforcement, sadeghi2023social}. These approaches not only naturally enforce resolving unsatisfied demands (as deprivation occurs during lack of goods or services), but also promote equity and acceptable response times due to the superlinear growth of deprivation costs over time. An overview of research on deprivation costs in humanitarian logistics can be found in \cite{shao2020research}, which concludes that most works using deprivation costs are limited to deterministic settings. In stochastic problems, research mainly considers minimizing unmet demands, costs, or response times rather than deprivation costs \citep[except for the two-stage stochastic program of][]{moreno2018effective}. Stochastic dynamic multi-period allocation models that include deprivation costs still remain unexplored, yet it is an important problem class to study. Disaster areas are inherently dynamic and stochastic, and multi-period problems capture the intertemporal effects of deprivation during an operation \citep{perez2016inventory}. Deprivation costs offer the opportunity to balance the conflicting goals of minimizing operational costs and increasing beneficiaries' welfare. In summary, incorporating deprivation costs in objective functions bring stochastic dynamic problems in humanitarian operations much closer to practice, yet this setup is understudied as of now.

\subsection{Contributions}
\label{sec:lit:contributions}
Stochastic dynamic solution methods have shown superior results in humanitarian allocation problems compared to static variants, as they adapt their plans to the evolving nature of a disaster area. Various multi-period models have been proposed that include the uncertainty of demand and supply. Among these models, often unsatisfied demand is minimized besides logistics costs, whereas deprivation costs -- as explicit incorporation of human suffering in the objective -- have not yet been considered in stochastic dynamic allocation problems. Nonetheless, the consideration of deprivation is essential to focus on human suffering, especially in situations with scarcity of supplies. In humanitarian applications, learning-based heuristics have been developed, but these only consider limited action spaces and mostly deterministic problems in the humanitarian allocation context. Appendix \ref{app:literature} summarizes the literature related to our work, including works in the broader field of stochastic dynamic allocation problems, to which we refer in Section \ref{sec:solution} regarding our solution approach. We contribute to the literature by addressing the stochastic dynamic post-disaster inventory allocation problem (SDPDIAP) with uncertain supply and demand that aims to minimize transportation costs and deprivation costs. To this end, we design two reinforcement learning methods based on approximate dynamic programming that anticipate the evolving nature of the aftermath of a disaster, one designed for the highest performance and the other for excellent scalability. Additionally, we provide practical insights into dealing with scarcity of supplies during dynamic disaster response operations through the use of deprivation costs. Finally, we show how the deployment of cargo UAVs offers opportunities to further improve transport efficiency and reduce deprivation costs.

\section{Model Formulation}
\label{sec:model}
This section formalizes the model for our humanitarian allocation problem. Considering the multi-periodicity required to fully integrate deprivation costs during humanitarian operations \citep{holguin2013appropriate} and uncertainties and dynamics in a disaster area \citep{anaya2016models}, we describe the stochastic dynamic post-disaster inventory allocation problem (\mbox{SDPDIAP}) with uncertain supply and demand. The main challenge is to allocate supplies -- which are generally limited in disaster areas \citep{najafi2013multi} -- such that the operation is efficient in terms of minimizing transportation costs, but also effective in alleviating human suffering by minimizing deprivation costs. We first describe the problem in Section~\ref{subsec:problem} and then formulate the problem as a Markov decision process in Section~\ref{subsec:model}.

\subsection{Problem Description}
\label{subsec:problem}
We build upon the widely recognized relief supply chain as described by \cite{balcik2008last}, in which a central warehouse (CW) forms the main entry point for incoming supplies into the disaster area and allocates the supplies to local distribution centers in specific districts. From the local distribution centers, last-mile distribution is performed to deliver supplies to beneficiaries in the districts. We focus on the allocation problem of the CW towards the districts.

\begin{figure}[ht!]
  \centering
  \includegraphics[width=0.5\linewidth]{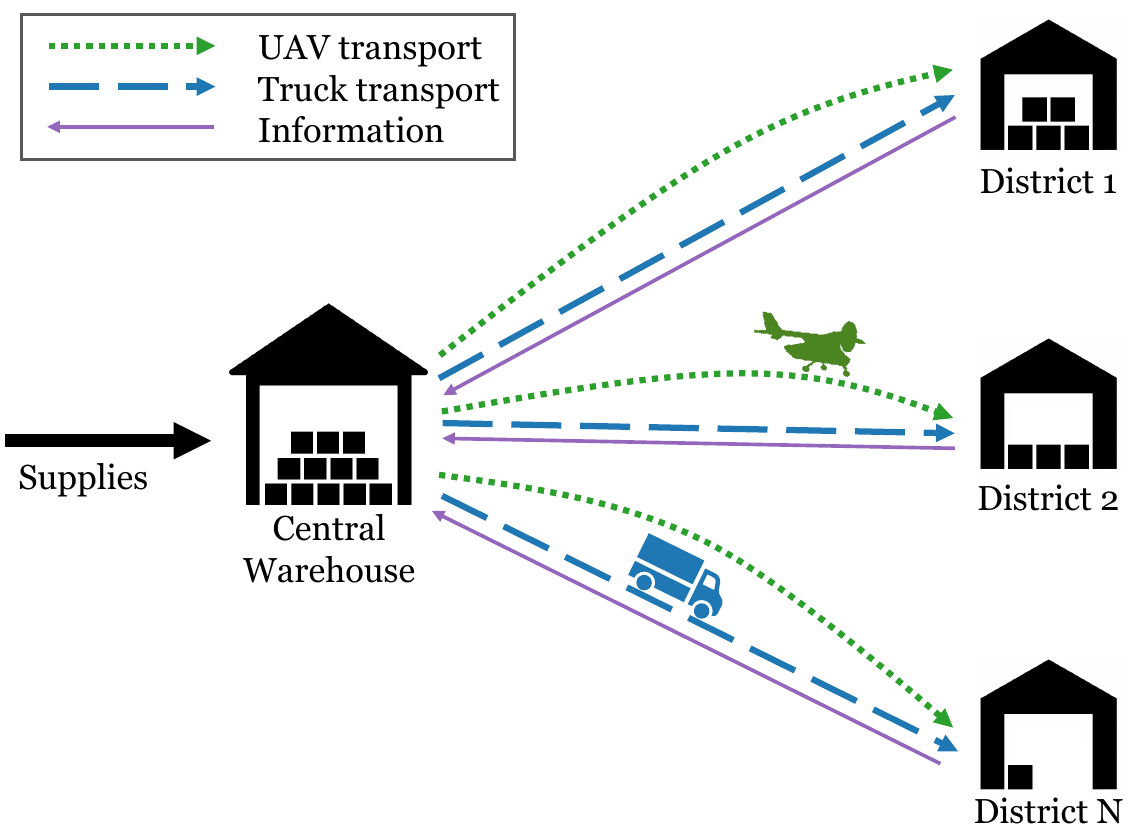}
  \caption{Illustration of the relief allocation problem. Supplies at the CW are allocated to districts with trucks or UAVs. Each district informs the CW regarding inventory, demand, and deprivation.}
  \label{fig:relief_chain}
\end{figure}

After the event of a disaster, relief supplies arrive at the CW over time according to a stochastic process. Allocation decisions for these supplies are made by the humanitarian agency operating the CW at decision epochs $t \in \set{T} = \{0, 1, ..., T\}$, e.g., every six hours. The relief supplies arriving throughout $[t-1, t)$ are captured by $\hat{s}^{\text{CW}}_t$ (realizations of random variables are indicated by the hat operator). The total available supply inventory $I^{\text{CW}}_t$ can be distributed to any district $n \in \set{N} = \{1, 2, ..., N\}$ by means of transportation mode $k \in \set{K} = \{1, 2, ..., K\}$ (more specifically, trucks or UAVs), but can also be kept in inventory for future decision epochs. Each transportation mode $k$ has an associated capacity $q_k$ and fixed transportation costs per district $n$ for one vehicle, denoted by $c_{nk}$, assuming that the costs of sending a vehicle to a district is independent of its load. We allow multiple vehicles to visit a district at a single decision epoch. Allocated supplies $x_{tnk}$ replenish the district inventory $I_{tn}$, which are consumed throughout $[t-1, t)$ according to the random district demand $\hat{d}_{tn}$. At time $t$, the realizations of demand and inventory positions become known by information updates from aid workers in each district, who also provide estimates for the demand in the upcoming period, denoted by $d_{t,t+1,n}$. Figure \ref{fig:relief_chain} illustrates the problem of allocating supplies utilizing trucks and UAVs to different districts and the information flows back to the CW.

An essential component of the SDPDIAP is the integration of human suffering in the problem definition. We follow the concept of deprivation costs as defined by \cite{holguin2013appropriate} and also implemented by \cite{perez2016inventory} in a deterministic variant of the post-disaster inventory allocation problem. In a multi-period setting, deprivation cost functions are often discretized in fixed time steps that preserve the superlinear nature of the underlying continuous function \citep{rivera2016dynamic}. As we consider discrete decision epochs $t \in \set{T} = \{0, 1, ..., T\}$, we also discretize the nonlinear deprivation costs for each epoch. We apply the same deprivation cost function $\gamma(\delta)$ as function of the deprivation time $\delta$ as \cite{perez2016inventory}:
\begin{gather}
    \gamma(\delta) = e^{0.065\delta} - 1. \label{eq:depr_cost_func}
\end{gather}

Equation (\ref{eq:depr_cost_func}) represents the deprivation costs per individual per hour without the required supplies. Instead of hourly costs, we consider the costs of a period between decision epochs, where the length of a period is incorporated into the deprivation time $\delta_{tn}$ at epoch $t$ in district $n$. Once inventory in a district is depleted by the demand, people in the district will suffer during the deprivation time $\delta_{tn} \in \mathbb{Z}_{\geqslant 0}$ and deprivation costs arise according to the deprivation cost function $\gamma(\delta_{tn})$. The number of individuals for which the demand is unsatisfied is denoted by the shortage $h_{tn}$. If a fraction of the demand is satisfied, the shortage is corrected for that fraction of the demand throughout a period. The deprivation costs per individual in a district are multiplied by the shortage $h_{tn}$ to obtain the district deprivation costs. The increase of deprivation costs between consecutive time steps is referred to as the marginal deprivation cost \citep{ismail2021possibilistic}. We express the marginal deprivation cost between two decision epochs by $g(\delta_{tn}) = \max\{ 0, \gamma(\delta_{tn}) - \gamma(\delta_{tn} - 1) \}$. The deprivation time and resulting deprivation costs will drop to zero when new supplies arrive that are sufficient to fulfill all demand at point $t$, and only start increasing again when the inventory is depleted \citep{holguin2013appropriate}. When the demand is not sufficient, deprivation will only be accounted for the shortage, but the deprivation time is not reset to zero. Demand is non-additive (i.e., no back-orders) as one individual could not, e.g., eat or drink for three days when supplies finally arrive.

The objective is to minimize the transportation costs for all allocations and the total deprivation costs over all districts $n \in \set{N}$ and decision epochs $t \in \set{T}$. In this way, allocation decisions will be made according to the needs of each district, using a suitable transportation mode.

\subsection{Sequential Decision Process}
\label{subsec:model}
This section describes the stochastic dynamic post-disaster inventory allocation problem as a Markov decision process (MDP). For this, we use the framework proposed by \cite{powell2019unified} for sequential decision-making. Subsequently, we present the state variables, decision variables, exogenous information, transition function, and the objective function and optimal policy.

\paragraph{State Variables.}
The state $S_t$ at decision epoch $t \in \set{T}$ is defined by the tuple:
\begin{gather}
    S_t = \left(I^{\text{CW}}_t, \left( I_{tn}, h_{tn}, \delta_{tn}, d_{t,t+1,n} \right)_{n \in \set{N}} \right). \label{eq:state_var}
\end{gather}

The state contains all dynamic information that is necessary to compute the decision, the transitions, and the rewards. Specifically, $I^{\text{CW}}_t$ represents the inventory position at the central warehouse at $t$, and for each district $n \in \set{N}$, we have the inventory position $I_{tn}$, current shortage $h_{tn}$, deprivation time $\delta_{tn}$, and expected demand $d_{t,t+1,n}$ for the upcoming period $[t, t+1)$.

\paragraph{Decision Variables.}
At each decision epoch $t$, supplies at the central warehouse $I^{\text{CW}}_t$ can be allocated to districts $n$ with transportation modes $k$. Matrix $\mathbf{x}_t = \left[x_{tnk}\right]_{n \in \set{N}, k \in \set{K}}$ denotes the decision, describing the number of supplies allocated to each district $n$ by each transportation mode $k$ at epoch $t$. The sum of allocations is constrained by the number of supplies available in the CW. Relief supplies are not necessarily all allocated to districts, but can also be kept in inventory at the CW for future allocations. The set of feasible decisions in a state is described by:
\begin{gather}
    \set{X}(S_t) = \left\{ \mathbf{x}_t \in \mathbb{N}^{|\set{K}||\set{N}|} : \sum_{k \in \set{K}}\sum_{n \in \set{N}} x_{tnk} \leqslant I^{\text{CW}}_t   \right\}.
\end{gather}

\paragraph{Exogenous Information.}
Between subsequent decision epochs $t$ and $t+1$, the disaster area evolves and new information arrives. Realizations of random variables in the new information are indicated by the hat operator. The exogenous information process in our problem concerns the arrival of supplies in the central warehouse $\hat{s}^{\text{CW}}_{t+1}$, and for each district $n \in \set{N}$ the realization of demands $\hat{d}_{t+1,n}$ during $[t, t+1)$ and estimates for the demands $d_{t+1,t+2,n}$ in the upcoming period $[t+1, t+2)$, which are communicated by aid workers in the districts. Let $\omega \in \Omega$ describe a sample path $W_1, W_2,...,W_T$, with a particular realization of the exogenous information $W_{t+1}$ given by:
\begin{gather}
    W_{t+1} = \left( \hat{s}^{\text{CW}}_{t+1}, (\hat{d}_{t+1,n}, d_{t+1,t+2,n})_{n \in \set{N}}    \right).
\end{gather}

\paragraph{Transition Function.}
To model the system forward from $t$ to $t+1$, we define a transition function to the next state by: $S_{t+1} = S^M(S_t, \mathbf{x}_t, W_{t+1})$, incorporating the decision and new information. Without loss of generality, allocations from the CW at time $t$ are assumed to arrive in the same period at the districts (i.e., immediate increase of district inventory). Deprivation times increase by one if the inventories including allocations are not sufficient to fulfill the period demand and are otherwise reset to zero. The transition function updates the state such that:
\begin{align}
    I^{\text{CW}}_{t+1} & = I^{\text{CW}}_t - \sum_{k \in \set{K}}\sum_{n \in \set{N}} x_{tnk} + \hat{s}^{\text{CW}}_{t+1}, \label{eq:trans-Icw} \\
    I_{t+1,n} & = \max \left\{ 0, \enspace I_{tn} + \sum_{k \in \set{K}} x_{tnk} - \hat{d}_{t+1,n} \right\}, & \forall n \in \set{N}, \label{eq:trans-I} \\
    h_{t+1,n} & = \max \left\{ 0, \enspace \hat{d}_{t+1,n} - \enspace I_{tn} + \sum_{k \in \set{K}} x_{tnk} \right\}, & \forall n \in \set{N}, \label{eq:trans-alpha} \\
    \delta_{t+1,n} & = \begin{cases}
                        \delta_{tn} + 1 & \text{if } I_{tn} + \sum_{k \in \set{K}} x_{tnk} \leqslant \hat{d}_{t+1,n}, \\
                        0 & \text{otherwise,}
                        \end{cases} & \forall n \in \set{N}, \label{eq:trans-delta} \\
    d_{t+1,t+2,n} & = d_{t+1,t+2,n}, & \forall n \in \set{N}, \label{eq:trans-dem}
\end{align}

\noindent where the inventory position of the central warehouse is updated according to the allocation decisions and the arrival of new supplies in Equation (\ref{eq:trans-Icw}). The inventory positions at the districts are increased by the allocation decisions and decreased by the realized demand in the period or set to zero if the demand exceeds the inventory position plus the allocation decisions in Equations (\ref{eq:trans-I}). On the contrary, the shortage becomes positive if the demand exceeds the inventory plus the allocations in Equations (\ref{eq:trans-alpha}), or zero otherwise. The deprivation times are increased by one if the post-decision inventory is not sufficient to fulfill the demand and otherwise reset to 0 in Equations (\ref{eq:trans-delta}). Lastly, the new demand estimate is taken from the exogenous information in Equations (\ref{eq:trans-dem}).

\paragraph{Objective Function and Optimal Policy.}
The goal is to minimize the total deprivation costs and transportation costs over time over all districts. The costs at decision epoch $t$ are given by:
\begin{gather}
    C(S_t, \mathbf{x}_t) = \sum_{n \in \set{N}} g(\delta_{tn}) h_{tn} + \sum_{k \in \set{K}}\sum_{n \in \set{N}} \left \lceil \frac{x_{tnk}}{q_k} \right \rceil c_{nk}, \label{eq:costfunction}
\end{gather}

\noindent where the first part captures the deprivation costs $g(\delta)$ based on the deprivation time $\delta_{tn}$ and shortages $h_{tn}$ in the districts, and the second part captures the transportation costs for the required number of vehicles to each district, where the required number of vehicles is computed by the allocated supplies $x_{tnk}$ and the vehicle capacity $q_k$, multiplied by the costs $c_{nk}$ per vehicle per district. In the final state $S_T$ no decisions are made, but the terminal deprivation costs over the last period $[T-1, T)$ are calculated according to the cost function. 

Let $\pi$ be a policy in the set of policies $\Pi$, and let $X_t^\pi(S_t)$ be the decision function of $\pi$ that determines decision $\mathbf{x}_t \in \set{X}(S_t)$ for each state $S_t \in \set{S}$. We aim to find the optimal policy $\pi^*$ that minimizes the expected costs $C_t$ over all decision epochs $t \in \set{T}$ starting from the initial state $S_0$, given by the objective function:
\begin{gather}
    \min_{\pi \in \Pi} \mathbb{E} \left\{ \sum_{t=0}^T C\left(S_t, X_t^\pi(S_t)\right) \mid S_0 \right\}. \label{eq:optimalpolicy}
\end{gather}

Although computationally intractable, we could formally compute the optimal policy by solving the well-known Bellman equations \citep{bellman1957markovian}. To obtain a decision $\mathbf{x}_t$, we consider the direct costs $C(S_t, \mathbf{x}_t)$ and the expected future costs of the next state $S_{t+1}$ to account for the downstream effect of decisions. Let $V_t(S_t)$ denote the expected future costs of state $S_t$. The optimal decisions $\mathbf{x}^*_t$ at $t \in \set{T}$ can be found by:
\begin{gather}
    \mathbf{x}^*_t = \underset{\mathbf{x}_t \in \set{X}(S_t)}{\mathrm{argmin}} \left( C(S_t, \mathbf{x}_t) + \mathbb{E} \left\{ V_{t+1}(S_{t+1}) \mid S_t, \mathbf{x}_t \right\} \right),
\end{gather}

\noindent where $S_{t+1} = S^M(S_t, \mathbf{x}_t, W_{t+1})$. In the next section, we will describe our solution methods, which aim to find such a policy within reasonable time.

\section{Solution Methods}
\label{sec:solution}
In principle, the optimal policy can be obtained by solving Equation (\ref{eq:optimalpolicy}), however, this is computationally intractable in our problem setting, as we face the three curses of dimensionality \citep{powell2011approximate}. First, the outcome space is large, determined by the random information regarding supply and demand, and it increases exponentially with a higher number of districts. Second, we cannot learn the values of each individual state as the inventory positions and demands within a state can take many values. Third, the decision space is vast due to the combinatorial nature of the supply allocation across multiple districts with different transportation modes. We adopt the Approximate Dynamic Programming (ADP) framework of \cite{powell2011approximate} to address these problems and extend this with a decomposition approach and a neural network approach. ADP is a simulation-based solution method that aims to learn a Value Function Approximation (VFA) that gives the estimated future costs. Future costs are defined by the value of being in the next state $V_{t+1}(S_{t+1})$. The problem is solved by iteratively simulating the MDP and learning the value of being in a certain state by observing the costs. In Section \ref{subsec:outcome}, we describe how we deal with the high-dimensional outcome space; in Section \ref{subsec:state} we discuss how we handle the large state space and input features for the value functions; in Section \ref{subsec:action} we describe how we solve the combinatorial action space, and the overall algorithmic outline is presented in Section \ref{subsec:algo_outline}.

\subsection{Outcome Space}
\label{subsec:outcome}
To avoid enumerating the outcome space, ADP relies on the concept of the post-decision state \citep{powell2011approximate}. The post-decision state $S^\text{x}_t$ is the state immediately after decision $\mathbf{x}_t$ and before new information $W_{t+1}$, with $S^\text{x}_t = S^{Mx}(S_t, \mathbf{x}_t)$ describing the deterministic transition:
\begin{gather}
    S^\text{x}_t = \left(I^{\text{CW, x}}_t, \left( I^\text{x}_{tn}, h_{tn}, \delta_{tn}, d_{t,t+1,n} \right)_{n \in \set{N}} \right), \label{eq:pd_state}
\end{gather}

\noindent where $I^{\text{CW, x}}_t = I^{\text{CW}}_t - \sum_{k \in \set{K}}\sum_{n \in \set{N}} x_{tnk}$ and $I^\text{x}_{tn} = I_{tn} + \sum_{k \in \set{K}} x_{tnk} \enspace \forall n \in \set{N}$. From the post-decision state $S^\text{x}_t$ we can compute the next state $S_{t+1}$ with transition $S_{t+1} = S^{MW}(S^\text{x}_t, W_{t+1})$, incorporating exogenous information $W_{t+1}$. By computing the future costs for the post-decision state $S^\text{x}_t$ instead of the next state $S_{t+1}$, we can use the deterministic value function $V_t(S^\text{x}_t)$, instead of the expectation $\mathbb{E} \left\{ V_{t+1}(S_{t+1}) \mid S_t, \mathbf{x}_t \right\}$ over the outcome space of the future costs, avoiding enumeration over all possible supply- and demand realizations.

\subsection{State Space}
\label{subsec:state}
The state space is large, due to the many values that inventory positions and demands can take. For moderately sized problems, we must already perform an intractable number of computations to visit each post-decision state at least once to estimate the future costs. Furthermore, the state space grows significantly with each additional district to serve, increasing the state dimension with four variables. The general approach to deal with large state spaces is to substitute the true value function $V_t(S^\text{x}_t)$ with a value function approximation (VFA) $\overline{V}_t(S^\text{x}_t)$. The benefit of VFA is that it avoids full enumeration of the state space to obtain future costs of each state. We design two alternative VFA approaches: (i) a decomposed linear value function approximation (DL-VFA) aimed at scalability, which decomposes the state in district states and learns separate linear value functions for each district at each decision epoch, and (ii) a neural network VFA (NN-VFA) in which one neural network approximates all state values, designed to better capture nonlinearities. We briefly outline both methods.

\paragraph{VFA Methods.} Instead of approximating the future costs of a state with a single value function, when the problem structure allows one might decompose it into multiple approximations. For example, \cite{bouzaiene2016single} learn multiple VFAs per locomotive in a locomotive allocation problem, \cite{zhang2022online} introduce a value function for each parking area in a parking allocation problem, and \cite{beirigo2022learning} consider an autonomous mobility-on-demand system with different functions for owned- and hired vehicles. The advantages of decomposing the value function are improvement in the scalability of solution methods, reduction of complexity in the learning process, and clearer insights into the contribution per subfunction. With DL-VFA, we decompose the state $S^\text{x}_t$ into sub-states based on the districts $S^\text{x}_{tn} n \in \set{N}$ and approximate the district costs $\overline{V}_{tn}(S^\text{x}_{tn})$. This decomposition is possible because the state contains clear district-related variables and the cost function (\ref{eq:costfunction}) consists solely of district-related costs. The overall value function is obtained by summing the district value functions. The benefit of this decomposition is that the complexity of learning the cost structure of a single district is significantly lower than learning the future costs of all districts together, and with VFAs per district, allocation decisions could be explained by the function weights per district. As our problem has a finite horizon $T$ and is non-stationary (i.e., supply, demand, and expected value of demands $d_{t,t+1,n}$ may vary over time), future costs are time-dependent and computed for each epoch $t \in \set{T}$. Hence, we will learn value functions not only per district but also per decision epoch $t \in \set{T}$. With DL-VFA, the value functions are approximated by a linear combination of features based on the post-decision states, denoted by $\Phi(S^\text{x}_t)$ and described below, explaining the cost structure within a district. To accurately represent the future costs of districts, we learn the weights $\Theta_{tn}$ for each district and each epoch related to the feature vector $\Phi(S^\text{x}_{tn})$ derived from district state $S^\text{x}_{tn}$. We obtain the DL-VFA by:
\begin{gather}
    V_t(S^\text{x}_t) \approx \overline{V}^{\text{linear}}_t(S^\text{x}_t) = \sum_{n \in \set{N}} \overline{V}^{\text{linear}}_{tn}(S^\text{x}_{tn}) = \sum_{n \in \set{N}} \Theta_{tn}\Phi^{\text{linear}}(S^\text{x}_{tn}), \enspace \forall \enspace t \in \set{T}.
\end{gather}

In our objective function, we deal with the nonlinearity of deprivation costs. To include their possible impact on future costs for an improved approximation of the value function, we propose a neural network as a nonlinear VFA that approximates the future costs of a full state $S^\text{x}_t$ instead of decomposed district states $S^\text{x}_{tn}$, at the expense of computational efficiency and explainability of DL-VFA. For this neural network VFA (NN-VFA), denoted by $\mathbb{F}(S^\text{x}_t)$, we train the parameter matrix $\Theta^{\text{neural}}$ (i.e., weights in the neural network) to approximate the future costs based on the feature vector $\Phi^{\text{neural}}(S^\text{x}_t)$. We obtain the NN-VFA by:
\begin{gather}
    V(S^\text{x}_t) \approx \overline{V}^{\text{neural}}(S^\text{x}_t) = \mathbb{F}(\Phi^{\text{neural}}(S^\text{x}_t) | \Theta^{\text{neural}}).
\end{gather}

\paragraph{Feature Design.} The feature vectors $\Phi^{\text{linear}}$ and $\Phi^{\text{neural}}$, based on the post-decision (district) states, are used as input for the VFAs to provide information based on the post-decision states $S^\text{x}_{tn}$ and $S^\text{x}_t$, and aim to capture the problem's cost structure. The main feature we define is the expected deprivation cost $\overline{G}^{\text{x}}_{t, t+1, n}$. This is an important predictor, because deprivation is only reset to zero if the post-decision inventory at $t$ covers the demand during $[t, t+1)$. To obtain the expected deprivation costs for the next period, we first define the expected shortage $\overline{h}_{t, t+1, n}$ as the expected number of individuals that will suffer from a lack of supplies. The expected deprivation costs are the product of the expected shortage $\overline{h}_{t, t+1, n}$ and the individual deprivation in the next period $g(\delta_{tn} + 1)$. As the demand -- and thus the demand coverage -- is uncertain, we take a cautious view of the expected shortage. Hence, we add two standard deviations to the expected demand and compute the expected shortage as:
\begin{gather}
    \overline{h}_{t, t+1, n} = \max\left\{0, \enspace d_{t, t+1, n} + 2 {\sigma_d}_{t, t+1, n} - I^{\text{x}}_{tn} \right\}.
\end{gather}

In this way, the probability of an actual shortage is low when the post-decision inventory is increased such that the expected shortage $\overline{h}_{t, t+1, n}$ is eliminated, resulting in a high probability that deprivation is reset to zero. The expected deprivation costs only function as a feature for the VFA, the actual deprivation costs will be calculated based on the realized shortages in districts. The expected deprivation costs are obtained by:
\begin{gather}
    \overline{G}^{\text{x}}_{t, t+1, n} = g(\delta_{tn} + 1) \cdot \overline{h}_{t, t+1, n}.
\end{gather}

Other important predictors of the spatio-temporal impact on both transportation costs and deprivation costs are the inventory positions  $I^\text{x}_{tn}$ and deprivation times $\delta_{tn}$ in districts, which we therefore also include in the feature vectors. If the inventory in a particular district is high, both transportation and deprivation are not expected to occur in future periods (i.e., no allocations are required). If the inventory is zero, both transportation costs (to replenish the district) and deprivation costs (until replenishment) are expected. Furthermore, deprivation time will increase due to lack of supplies, increasing the deprivation costs over time. As the neural network of NN-VFA is a single value function for all decision epochs and the complete state, the time index $t$ and central warehouse inventory $I^{\text{CW,x}}_t$ are provided as features for the neural network. In summary, the feature vectors for the linear value functions and the neural network are defined as:
\begin{align}
    \Phi^{\text{linear}}(S^\text{x}_{tn}) & = \left( I^\text{x}_{tn}, \delta_{tn}, \overline{G}^{\text{x}}_{t, t+1, n}\right), \enspace & \forall \enspace n \in \set{N}, \enspace t \in \set{T} \label{eq:lin_feat}, \\
    \Phi^{\text{neural}}(S^\text{x}_t) & = \left( t, I^{\text{CW,x}}_t,  I^\text{x}_{tn}, \delta_{tn}, \overline{G}^{\text{x}}_{t, t+1, n} \right)_{n \in \set{N}}, \enspace & \forall \enspace t \in \set{T} \label{eq:neural_feat}.
\end{align}

\subsection{Combinatorial Action Space}
\label{subsec:action}
At each decision epoch, we decide upon the combination of the supply quantities for (i) each district and (ii) each transportation mode. Such combinatorial action spaces are common in transportation, but they pose a challenge for standard reinforcement learning algorithms \citep{dulac2021challenges} that require full enumeration of the action space. Furthermore, standard algorithms can often not handle constrained action spaces, resorting to masking to filter out infeasible actions.

To efficiently select actions in the large constrained action space in DL-VFA and NN-VFA, we formulate the single-stage decision problem as a mixed-integer program (MIP) and solve this using a general-purpose solver. The MIP solves the approximate Bellman equation and handles the constraints, such as the limited supplies in the central warehouse. In this case, the value function approximations (VFAs) are directly integrated into the MIP. This immensely increases the potential action spaces we can handle, as the single-stage decision problem can be solved in a fraction of the time required for enumeration, while preserving optimality guarantees in selecting the best action under the prevailing policy. This approach, using a single (piecewise) linear VFA in each stage, is applied in stochastic dynamic problems such as intermodal transport \citep{rivera2016dynamic, heinold2022primal}, ridesharing systems \citep{lei2019path, beirigo2022learning}, and delivery dispatching problems \citep{klapp2018one, heeswijk2019delivery}. For the DL-VFA, our methodological contribution lies in the decomposition of the state and accompanying VFA. The decomposition yields one VFA for each district and each epoch. The VFAs for all districts in one epoch are then integrated into the single-stage decision problem. By solving this problem, the optimal combination of allocations is made given the district-based value functions. 

The objective incorporates the direct costs $C(S_t, \mathbf{x}_t)$ from Equation (\ref{eq:costfunction}) and the approximate future costs $\overline{V}_t(S^\text{x}_t)$, based on the post-decision state features. The objective of DL-VFA at each decision epoch includes the direct deprivation and transportation costs and the approximated future costs per district:
\begin{align}
    \min_{x_{tnk} \in \set{X}} \sum_{n \in \set{N}} \left( g(\delta_{tn}) \hat{d}_{tn} + \sum_{k \in \set{K}} \left \lceil \frac{x_{tnk}}{q_k} \right \rceil c_{nk} + \Theta_{tn}\Phi^{\text{linear}}(S^\text{x}_{tn})  \right) \label{eq:lin_obj},
\end{align}

\noindent where $\Theta_{tn}$ are the trainable weights and $\Phi^{\text{linear}}$ the input features according to Equation (\ref{eq:lin_feat}). NN-VFA includes the approximated future costs as output of the neural network in the objective:
\begin{align}
    \min_{x_{tnk} \in \set{X}} \sum_{n \in \set{N}} \left( g(\delta_{tn}) \hat{d}_{tn} + \sum_{k \in \set{K}} \left \lceil \frac{x_{tnk}}{q_k} \right \rceil c_{nk} \right) + \mathbb{F}(\Phi^{\text{neural}}(S^\text{x}_t) | \Theta^{\text{neural}}) \label{eq:nn_obj},
\end{align}

\noindent where $\mathbb{F}$ is the neural network with trainable weights $\Theta^{\text{neural}}$ and input features $\Phi^{\text{neural}}$ according to Equation (\ref{eq:neural_feat}). Both models are subject to the following constraints:
\begin{align}
    I^{\text{CW,x}}_t & = I^{\text{CW}}_{t-1} - \sum_{k \in \set{K}}\sum_{n \in \set{N}} x_{tnk}, \label{eq:lin_cw} \\
    I^\text{x}_{tn} & = I_{t-1,n} + \sum_{k \in \set{K}} x_{tnk}, & \forall \enspace n \in \set{N}, \label{eq:lin_inv} \\
    \overline{G}^{\text{x}}_{tn} & = g(\delta_{tn} + 1) \overline{h}_{t, t+1, n}, & \forall \enspace  n \in \set{N}, \label{eq:lin_depr} \\
    \overline{h}_{t, t+1, n} & \geqslant d_{t, t+1, n} + \xi_{t, t+1, n} - I_{tn}, & \forall \enspace n \in \set{N}, \label{eq:lin_exp_short_a} \\
    \overline{h}_{t, t+1, n} & \leqslant d_{t, t+1, n} + \xi_{t, t+1, n} - I_{tn} -M^{\text{h}}_{n} \cdot (1 - z^{\text{h}}_n), & \forall \enspace n \in \set{N}, \label{eq:lin_exp_short_b} \\
    \overline{h}_{t, t+1, n} & \leqslant M^{\text{h}}_{n} \cdot z^{\text{h}}_n, & \forall \enspace n \in \set{N}, \label{eq:lin_exp_short_c} \\
    I^{\text{CW}}_t, \enspace & x_{tnk}, \enspace I_{tn}, \enspace h_{tn}, \enspace \overline{G}^{\text{x}}_{tn}, \enspace \overline{h}_{t, t+1, n} \geqslant 0, & \forall \enspace t \in \set{T}, \enspace n \in \set{N}, \enspace k \in \set{K}, \label{eq:lin_sets} \\
    z^{\text{h}}_n & \in \{ 0, 1 \} \label{eq:lin_bin}.
\end{align}

\noindent Equation (\ref{eq:lin_cw}) determines the post-decision inventory position at the central warehouse and Equations (\ref{eq:lin_inv}) compute the post-decision inventory positions in the districts. Equations (\ref{eq:lin_depr}) compute the expected deprivation costs in each district. Equations (\ref{eq:lin_exp_short_a})-(\ref{eq:lin_exp_short_c}) compute the expected shortage with big-M value $M^{\text{h}}_{n}$. Variable sets are defined by Equations (\ref{eq:lin_sets})-(\ref{eq:lin_bin}).

The model of DL-VFA is complete with these constraints as the linear value functions are directly integrated into the objective of the MIP in Equation (\ref{eq:lin_obj}). The objective function of NN-VFA, Equation (\ref{eq:nn_obj}), only contains the output of the neural network. The network itself needs to be integrated into the MIP as a set of constraints, following, e.g.,  \cite{delarue2020reinforcement} and \cite{van2020deep}. A neural network is characterized by $i \in \set{I} = \{ 1,..., I \}$ layers with $j \in \set{J}^i = \{ 1,..., J^i \}$ neurons, where the output of a neuron is denoted as $m^i_j$. Each neuron computes the dot product of the output vector of previous' layer neurons $\mathbf{m}^{i-1}$ with a vector of trainable weights $\theta^i$, adds a bias $\beta^i_j$, and computes the activation. We apply ReLU activation functions, where $ReLU(x) = \max\{x, 0\}$. The maximization statement of the ReLU activation may be written with big-M constraints where $M^-_{ij}$ and $M^+_{ij}$ denote the lower and upper bound of the function $\theta^i \cdot \mathbf{m}^{i-1} + \beta^i_j$ \citep{anderson2020strong}. Based on preliminary testing, we use a fully connected feed-forward neural network with two hidden layers to capture the complexity of our problem. The neural network value function $\mathbb{F}(\Phi^{\text{neural}}(S^\text{x}_t) | \Theta^{\text{neural}})$ with two hidden layers can be written using the following constraints for the single-stage decision problem: 
\begin{align}
    m^1_j & \geqslant \Theta^1 \cdot \Phi^{\text{neural}} + \beta^1_j, & \forall j \in \set{J}^1, \label{eq:nn_1a}\\
    m^1_j & \leqslant \Theta^1 \cdot \Phi^{\text{neural}} + \beta^1_j  - M^-_{1j} \cdot (1-z_{ij}), & \forall j \in \set{J}^1, \label{eq:nn_1b}\\
    m^1_j & \leqslant M^+_{1j} \cdot z_{ij}, & \forall j \in \set{J}^1, \label{eq:nn_1c}\\
    m^2_j & \geqslant \Theta^2 \cdot \mathbf{m}^{1} + \beta^2_j, & \forall j \in \set{J}^2, \label{eq:nn_2a}\\
    m^2_j & \leqslant \Theta^2 \cdot \mathbf{m}^{1} + \beta^2_j  - M^-_{2j} \cdot (1-z_{ij}), & \forall j \in \set{J}^2, \label{eq:nn_2b}\\
    m^2_j & \leqslant M^+_{2j} \cdot z_{ij}, & \forall j \in \set{J}^2, \label{eq:nn_2c}\\
    \mathbb{F} & = \Theta^3 \cdot \mathbf{m}^{2} + \beta^3_1, \label{eq:nn_out}\\
    z_{ij}     & \in \{0,1\}, & \forall \enspace j \in \set{J}^i, \enspace i \in \set{I}, \label{eq:nn_misc}
\end{align}

\noindent where Equations (\ref{eq:nn_1a})-(\ref{eq:nn_1c}) form the first hidden layer with the feature vector $\Phi^{\text{neural}}$ as input, Equations (\ref{eq:nn_2a})-(\ref{eq:nn_2c}) from the second hidden layer, Equation (\ref{eq:nn_out}) provides the output of the neural network, and Equation (\ref{eq:nn_misc}) defines the binary variable.

\subsection{Algorithmic Outline}
\label{subsec:algo_outline}
This section describes the algorithms to learn the value function approximations in DL-VFA and NN-VFA. Our approach is based on approximate dynamic programming (ADP), which iteratively samples the stochastic process to update value functions and improve future decisions \citep{powell2011approximate}. We describe this process for DL-VFA, detailed in Algorithm \ref{alg:DL-VFA}, and discuss the modifications for NN-VFA. 

In contrast to regular ADP, we provide a warm start with reasonable estimates for the initial value functions. Preliminary experiments showed that random initial decisions complicate convergence, as multiple value functions are included in the MIP and poor estimates preclude properly weighing the value functions against each other. To obtain initial estimates, we generate a large set of trajectories with post-decision states and costs using a warm-up heuristic and estimate initial value functions with linear regression based on this set. We refer to this set of trajectories as the initial experience buffer $\set{B}_0$. The warm-up heuristic randomly allocates 1-3 UAVs to districts where the deprivation time is higher than $Z_1 \sim \text{unif}\{1,3\}$ and thereafter selects a random district $n_{truck}$ to send a truck with all available supplies if the deprivation time is higher than $Z_3 \sim \text{unif}\{1,3\}$. The pseudo-code of the warm-up heuristic can be found in Appendix \ref{app:warm-up}.
\begin{figure}[!ht]
\centering
\begin{minipage}{1.0\linewidth}
\begin{algorithm}[H]
    \footnotesize
    \linespread{1}\selectfont
    \caption{Decomposed linear VFA algorithm}
    \label{alg:DL-VFA}
    \begin{algorithmic}
        \State \textbf{Step 0.} Set buffer size $B$ and discount factor $\lambda$.
        \StateLong{\textbf{Step 1.} Initialize exploration parameter $\epsilon^0$, update frequency $u$, smoothing factor $\alpha^0$, and buffer $\mathcal{B}^0$ for $B$ horizons to fill the experience buffer by running the warm-up heuristic.}
        \StateLong{\textbf{Step 2.} Compute initial weights $\Theta^{\text{linear},0}_{tn} \enspace \forall \enspace t \in \mathcal{T}, \enspace n \in \mathcal{N}$ with least squares regression based on $\mathcal{B}^0$.}
        \For{$r = 1, 2, \dots, R$}
            \State \textbf{Step 3.} Generate sample path $\omega^r$.
            \State \textbf{Step 4.} Set the initial state to $S_0$ using $W_0 \in \omega^r$.
            \For{$t = 0, 1,..., T $}
                \State \textbf{Step 5a.} Obtain decision $\mathbf{x}_t$:
                \If{$Z \sim \textnormal{unif}(0,1) < \epsilon^r$}
                    \State Get $\mathbf{x}_t$ by using step 3a from the warm-up heuristic.
                \Else
                    \StateLong{Get $\mathbf{x}_t$ by solving the optimization problem (\ref{eq:lin_obj}), (\ref{eq:lin_cw})-(\ref{eq:lin_bin}) using VFA weights $\Theta^{\text{linear},r-1}_{tn} \enspace \forall n \in \mathcal{N}$.}
                \EndIf
                \StateLong{\textbf{Step 5b.} Compute and store post-decision district states $S^\text{x}_{tn} = S^{Mx}_n(S_{tn}, \mathbf{x}_{tn})$ and district costs $C_n(S_{tn}, \mathbf{x}_{tn}) \enspace \forall n \in \mathcal{N}$ in buffer $\mathcal{B}^r$.}
                \StateLong{\textbf{Step 5c.} Obtain new information $W_{t+1} \in \omega^r$ and transition to the next state $S_{t+1} = S^{MW}(S^\text{x}_t, W_{t+1})$.}
            \EndFor
            \If{$r$ is multiple of $u$}
                \StateLong{\textbf{Step 6a.} Remove outliers from buffer $\mathcal{B}^r$ if direct costs in the sample path are higher than $Q3 + 1.5IQR$.}
                \StateLong{\textbf{Step 6b.} Recursively compute realized future district costs $\hat{V}_{tn}$ with algorithmic discount factor $\lambda \quad \forall \enspace t \in \mathcal{T}, \enspace n \in \mathcal{N}$.}
                \StateLong{\textbf{Step 6c.} Obtain weights $\hat{\Theta}^{\text{linear},r}_{tn} \enspace \forall \enspace t \in \mathcal{T}, \enspace n \in \mathcal{N}$ with least squares regression based on $\Phi^{\text{linear}}_{tn} \mathcal{B}^r$.}
                \StateLong{\textbf{Step 6d.} Update weights of all value functions with $\Theta^{\text{linear},r}_{tn} = (1 - alpha^r) \Theta^{\text{linear},r-1}_{tn} + \alpha^r \hat{\Theta}^{\text{linear},r}_{tn} \enspace \forall t \in \mathcal{T}, n \in \mathcal{N}$.}
                \State \textbf{Step 6e.} Update $\epsilon^r$ and $\alpha^r$.
            \EndIf
        \EndFor
        \State \textbf{Step 7. return} value functions $\left\{ \overline{V}^R_{tn} \enspace \forall \enspace t \in \mathcal{T}, \enspace n \in \mathcal{N} \right\}$.
    \end{algorithmic}
\end{algorithm}
\end{minipage}
\end{figure}

After obtaining the initial value functions, $R$ training episodes of one horizon are performed during which the value functions will be updated. With every episode $r$, a new sample path $\omega^r$ is generated, and a full horizon is simulated. At every decision epoch $t \in \{0, 1, ..., T\}$, we obtain a decision by solving the single-stage decision problem. To acquire a diverse set of decisions to learn from and to avoid getting stuck in local optima during training, with probability $\epsilon^r$ a decision is obtained by using the decision-making approach of the warm-up heuristic, and $\epsilon^r$ decays after each value function update. After obtaining the decision, the post-decision state and direct costs are computed and stored in the experience buffer, replacing the oldest experience in a first-in-first-out manner. New observations $W_{t+1} \in \omega^r$ are used to transition to the next state. This process is repeated for $R$ episodes. Every $u$ episodes (e.g., every 10 episodes), weights are updated based on the buffer $\set{B}^r$. Our approach differs from regular ADP, in which the value function is updated after every epoch or episode based on the last observations, by updating based on batches containing a set of observations (i.e., the experience buffer). We obtain a steady learning process by gradually refreshing the experience buffer $\set{B}$, estimating intermediate weights using all observations in $\set{B}$, and gradually updating weights of the value functions by the decaying factor $\alpha^r$. Before the weight update, we remove outliers from the experience buffer to prevent large disruptions in the learning process. A sample path is removed if direct costs in the sample path are 1.5 times the interquartile range above the third quartile. After removing outliers, the realized future costs of each post-decision (district) state are computed recursively. Intermediate weights for value functions are computed with least squares regression using experience buffer $\set{B}^r$. The new value function weights are obtained by smoothing the old weights with the intermediate weights using smoothing factor $\alpha^r$, which follows the approach of recursive estimation using multiple observations described by \cite{powell2022reinforcement}. Afterward, $\epsilon^r$ and $\alpha^r$ are decreased according to a decay factor. Training NN-VFA is based on similar principles as training DL-VFA. Instead of multiple VFAs, one neural network is trained, also initialized on buffer $\set{B}^0$. The main difference lies in updating the value function weights (Steps 6c and 6d in Algorithm \ref{alg:DL-VFA}); these steps are replaced by backpropagation through the neural network, also using the buffer $\set{B}^r$. NN-VFA is described in Appendix \ref{app:NN-VFA}. In the next section, we describe how we set up experiments to test the performance of DL-VFA and NN-VFA.

\section{Experimental Setup}
\label{sec:experiments}
This section describes the experimental setup. First, we describe the settings of theoretical instances in Section \ref{subsec:test} and then the case study of the Nepal 2015 Earthquake in Section \ref{subsec:Nepal}. Finally, we describe different benchmark methods in Section \ref{subsec:benchmark}.

\subsection{Theoretical Instances}
\label{subsec:test}
Methods are tested using several theoretical instances, with 30 periods of 6 hours (about one week), one central warehouse, and 1 to 6 districts portraying different levels of complexity. Districts have distinct demand levels and transportation costs, and both supply and demand follow normal distributions with a coefficient of variation (CoV) of 0.2. As transportation modes, we consider (i) trucks that carry 5000 units of supplies and (ii) UAVs that can load 200 units of supplies, inspired by Wings For Aid MiniFreighter \citep{bamsey2021dossier}. To accurately characterize UAVs, we collaborate in this study with the Wings For Aid foundation, which develops a UAV system that delivers humanitarian goods to people affected by natural disasters and man-made crises. The specific demand levels and transportation costs for different numbers of districts are listed in Appendix \ref{app:instances} in Table \ref{tab:app:details_districts}.

For in-depth analysis, we evaluate five scenarios using the 3-district instance. In the first scenario, (i) only trucks are used for supply allocation, serving as a baseline to quantify the added value of UAVs later on. In the second and third scenario, different levels of uncertainty are evaluated: (ii) a low-uncertainty scenario with a CoV of 0.1, and (iii) a high-uncertainty scenario with a CoV of 0.3. The fourth and fifth scenario consider two distinct S-curved demand patterns: (iv) a logistic decline of demand starting at 150\% of the supply rate and ending at 50\% of the supply rate, and (v) a logistic growth of demand starting at 50\% of the supply rate and ending at 150\% of the supply rate at the end of the horizon. These scenarios either portray a surge in demand at the start of the horizon or a situation in which a growing amount of demand arises over time. During training and solving of these instances, we assume the demand pattern is known, while in the exogenous information, only the expected next period demand is provided. In all scenarios, the average supply per period is assumed equal to the average demand per period over the whole horizon. At the start of the horizons, this will thus create a shortage in supplies for the decreasing demand pattern, and an oversupply for the increasing demand pattern. The specific demand levels and transportation costs for the different scenarios can be found in Appendix \ref{app:instances} in Table \ref{tab:app:demand_patterns_variants}.

\subsection{Nepal 2015 Earthquake Case}
\label{subsec:Nepal}
\begin{figure}[ht!]
    \centering
    \includegraphics[width=3.5in]{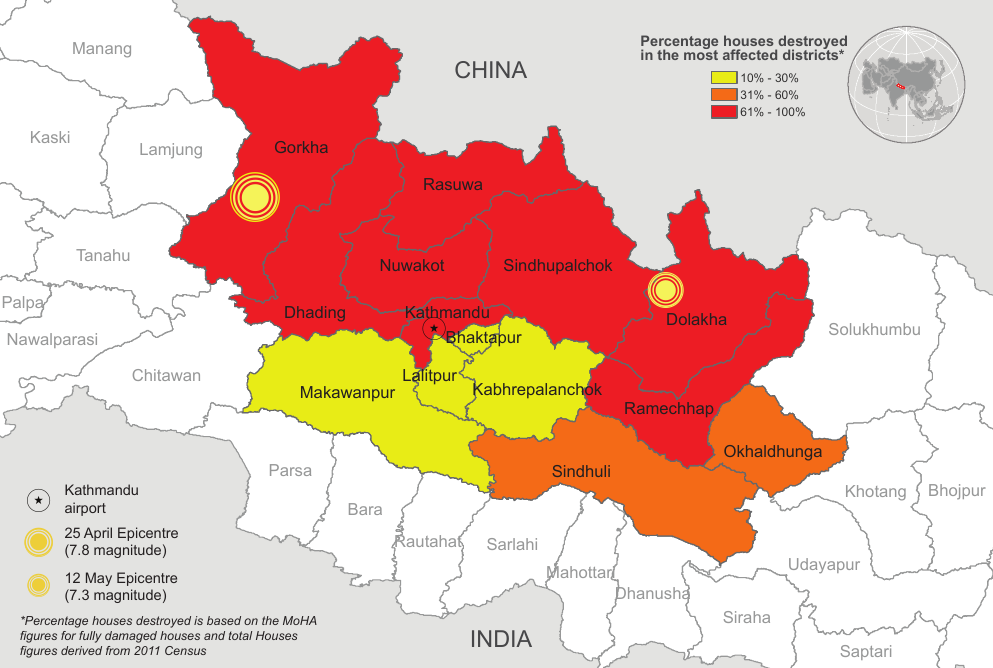}
    \caption{Disaster area overview, based on the OCHA Nepal Earthquake Humanitarian Snapshot of 24 July 2015.} \label{fig:Nepal_map}
\end{figure}

In April 2015, a 7.8M earthquake hit Nepal and another 7.3M earthquake two weeks later added to the catastrophe. Thirteen districts outside the capital city Kathmandu were heavily affected, with the mountainous surroundings of several districts hampering road transport. A main central warehouse was set up at the airport of Kathmandu, which received supplies from all over the world. The Nepal Red Cross Society coordinated the disaster response, and local offices in the thirteen districts set up distribution centers for the last-mile delivery of supplies in their districts. An overview of the disaster area and impact on houses is shown in Figure \ref{fig:Nepal_map}. We use data provided in \cite{johnsson2016tracking} regarding the thirteen districts. Truck costs are based on (i) the road distance from Kathmandu to the district capitals, (ii) the average geographical slope in the district, and (iii) the road density, in which the latter two indicate the level of difficulty navigating towards the districts. UAV costs are based on Euclidean distances. Demand is relative to the damaged houses, assuming that every 24 hours, 2\% of the damaged houses require supplies. The specific demand levels and transportation costs for the Nepal case are listed in Appendix \ref{app:instances} in Table \ref{tab:app:details_Nepal_instance}.

\subsection{Benchmark Methods}
\label{subsec:benchmark}
To analyze the performance of the proposed methods, we compare the results with a range of distinctive solution methods. As a lower bound, we solve a Perfect Information relaxation, and as comparative approaches, we compare the proposed methods with an exact re-optimization approach, a rule-based heuristic, and Proximal Policy Optimization (PPO), a state-of-the-art reinforcement learning method. We discuss these methods in the following paragraphs.

\paragraph{Perfect Information Bound.}
We consider the performance of an optimal policy under Perfect Information (PI), in which we presume all uncertainty to be resolved a priori. Instead of solving the problem one period at a time, the PI method solves all single-stage decision problems (with Equations (\ref{eq:lin_cw})-(\ref{eq:lin_bin})) at once, with the transition function (Equations (\ref{eq:trans-Icw})-(\ref{eq:trans-dem})) as constraints to tie the periods together. Whereas for other methods the sample path $\omega$ is revealed dynamically over time, for computing the PI bound we observe the full path a priori. Hence, this approach serves as a lower bound for the other methods. The objective follows the cost function, Equation (\ref{eq:costfunction}), and aims to minimize the costs over the complete horizon $\set{T}$:
\begin{gather}
    \min_{\mathbf{x} \in \set{X}} \sum_{t \in \set{T}} C(S_t,\mathbf{x}_t) = \min \sum_{t \in \set{T}} \left( \sum_{n \in \set{N}} g(\delta_{tn}) h_{tn} + \sum_{k \in \set{K}}\sum_{n \in \set{N}} \left \lceil \frac{x_{tnk}}{q_k} \right \rceil c_{nk} 
    \right). \label{eq:PI_obj}
\end{gather}

\paragraph{Re-optimization Approach.}
Each decision epoch, the re-optimization approach solves the same MIP model as the Perfect Information bound, but with expected values instead of sample realizations and only for the remaining time horizon $[t, \ldots, T]$. The model is updated every decision epoch by incorporating the previous decision and current inventory positions based on realized supply and demand. Thus, it offers exact solutions given the knowledge available at the current epoch. To accommodate demand uncertainty with this deterministic model, we account for an additional margin $2 {\sigma_d}_{t, t+1, n}$ for the next period demand, similar to the VFA approaches. Re-optimization follows the objective of the Perfect Information bound in Equation (\ref{eq:PI_obj}), but only for the remaining horizon $[t, \ldots, T]$.
 
\paragraph{Rule-based Heuristic.}
To mimic a practical situation without the use of advanced algorithms or solvers, we evaluate a rule-based heuristic. This heuristic is comparable to the warm-up agent of DL-VFA and NN-VFA, but without stochasticity in decision-making. For each decision epoch, if the deprivation time $\delta_n$ at district $n$ is 2, then sufficient UAVs are allocated to cover the expected demand. Thereafter, the district with the highest deprivation costs is selected and the UAV allocation plus all remaining goods in the CW are allocated to that district using trucks.

\paragraph{Proximal Policy Optimization.}
Proximal Policy Optimization (PPO) is a state-of-the-art reinforcement learning algorithm proposed by \cite{schulman2017proximal} that combines concepts of advantage actor-critic methods and trust region policy optimization. It has shown excellent performance in various problems, including several MuJoCo environments (robotic tasks with continuous actions), humanoid control tasks, and inventory management problems \citep{schulman2017proximal, vanvuchelen2020use}. PPO uses an actor-critic structure with an actor network $\mathbb{F}^a(S_t | \Theta^a)$ and a critic network $\mathbb{F}^c(S_t | \Theta^c)$. Networks are randomly initialized and then multiple Monte Carlo simulations are performed to fill an experience buffer used to update the policies. After each policy update, the experience buffer is emptied, as it is an on-policy method. The updates of both networks are clipped to provide relatively stable training behavior.

To adapt PPO to our problem, we let the actor network $\mathbb{F}^a$ output an action $\mathbf{x}^a_t$ directly (i.e., without solving the single-stage decision problem) with an action space of $1 + N K$ variables between $[0, 1]$. In the action space, one action is reserved for keeping supplies at the central warehouse, and the other actions concern each district-vehicle combination. As the actual action space varies based on the number of supplies in the central warehouse, we transform these $[0, 1]$ outputs of the actor network to actual supply allocations using a mapping function. The mapping function computes the ratio of each action and the sum of all actions and multiplies this with the available inventory in the central warehouse. To prevent small deliveries by UAVs, the amounts are rounded down to the nearest full UAV load, such that only full UAVs will be allocated to districts. In this way, the action $\mathbf{x}^a_t$ is transformed to the actual feasible action $\mathbf{x}_t$:
\begin{align}
    x_{tn,\text{UAV}} & = \left \lfloor \frac{I^{\text{CW}}_t x^a_{tnk}}{q_\text{UAV} \sum_{i \in \set{N}}\sum_{k \in \set{K}} x^a_{tik}} \right \rfloor q_\text{UAV}, & \forall \enspace n \in \set{N}, \\
    x_{tn,\text{truck}} & = \frac{I^{\text{CW}}_t x^a_{tnk}}{\sum_{i \in \set{N}}\sum_{k \in \set{K}} x^a_{tik}}, & \forall \enspace n \in \set{N}, \\
    \mathbf{x}_t & = (x_{tn,\text{UAV}}, \enspace x_{tn,\text{truck}})_{n \in \set{N}}.
\end{align}
\noindent Having described the experimental setup and benchmark methods, we proceed to the analysis of the experimental results.

\section{Numerical Results}
\label{sec:results}
This section presents the numerical results of DL-VFA, NN-VFA, and the benchmark policies. We analyze the results of different problem sizes in Section \ref{subsec:num_districts} and of different problem variants in Section \ref{subsec:num_variants}. Section \ref{subsec:num_Nepal} discusses the results of the Nepal case. We implemented all methods using Python 3.11 and Gurobi 10.0.2. Experiments were executed on a 32-core 2.6GHz Intel\raisebox{1ex}{\tiny{\textregistered}} Xeon\raisebox{1ex}{\tiny{\textregistered}} Platinum 8272CL machine with 64GB RAM. We cap training time at 4 hours for the learning-based methods, and solving time at 4 hours for the PI bound and at 15 minutes for each decision epoch of the re-optimization method, granting each solution method roughly equal computational budget. Specific parameter settings for each method are summarized in Appendix \ref{app:settings}.

\subsection{Results for Different Problem Sizes}
\label{subsec:num_districts}
Table \ref{tab:results_districts} shows detailed results for all methods in 1-6 districts. The best results are obtained by NN-VFA, outperforming the re-optimization benchmark by 7.3\% on average. Nonetheless, the computational effort is one of the highest, and with 6 districts -- the largest problem in this set -- the relative performance gain is only 2.2\%. DL-VFA outperforms the benchmark policies as well, improving on re-optimization by 6.2\%. DL-VFA is computationally efficient (almost 20 times faster than NN-VFA) and provides the best performance on the largest instance (7.8\% improvement on re-optimization), showcasing excellent scalability to larger problems. For all instances, the rule-based heuristic performs worse compared to DL-VFA, NN-VFA, and re-optimization. PPO performs reasonably well on small instances, but due to the combinatorially increasing action space performance drops below the rule-based heuristic performance in larger instances. Except for the 1-district instance, DL-VFA and NN-VFA match the maximum deprivation times and demand coverages of the PI bound, while obtaining costs of only 5.5\%-9.2\% higher than the bound. For UAV deployment, the PI bound spends around 37\% of the transportation costs on UAVs, while re-optimization spends 50\%, and the learning-based heuristics even 60\%. This shows that UAVs are deployed to adapt to supply and demand uncertainty, a consideration ignored by the PI bound.
\begin{table}[ht!]
\scriptsize
\centering
\caption{Results for each approach per instance size. The NN-VFA approach mostly outperforms all other methods, while both NN-VFA and DL-VFA outperform all benchmarks. The gap for re-optimization is the average gap over the 30 solutions in the horizon.}
\label{tab:results_districts}
\begin{adjustbox}{center}
\begin{tabular}{lrrrrrrrr}
\hline
\multirow{2}{*}{\textbf{1 district}} & \multicolumn{1}{r}{\multirow{2}{*}{Runtime (s)}} & \multicolumn{1}{r}{\multirow{2}{*}{Gap}} & Total costs                     & \multicolumn{1}{r}{Deprivation} & \multicolumn{1}{r}{UAV} & \multicolumn{1}{r}{Truck} & \multicolumn{1}{r}{Max. deprivation} & \multicolumn{1}{r}{Demand} \\
                  & \multicolumn{1}{r}{}         & \multicolumn{1}{r}{} & \multicolumn{1}{r}{(±std.dev.)} & costs       & costs & costs                   & time (h)            & coverage                   \\ \hline
PI-bound          &   11.99  & 0.00\%  & 3238 (±233)          & 1483        & 855   & 900  & 16     & 62\%   \\
Re-optimization   &   19.97  & 0.00\%  & 3876 (±120)          & 1656        & 1320  & 900  & 19     & 63\%   \\
Rule-based        &    0.01  &         & 4974 (±569)          & 924         & 0     & 4050 & 12     & 93\%   \\
PPO               & 3411.16  &         & 4034 (±165)          & 1109        & 2925  & 0    & 17     & 69\%   \\
DL-VFA            &  151.05  &         & \textbf{3536 (±335)} & 1361        & 1245  & 900  & 23     & 67\%   \\
NN-VFA            & 2625.25  &         & 3544 (±229)          & 1324        & 1320  & 900  & 22     & 67\%   \\
\hline
\multirow{2}{*}{\textbf{2 districts}} & \multicolumn{1}{r}{\multirow{2}{*}{Runtime (s)}} & \multicolumn{1}{r}{\multirow{2}{*}{Gap}} & Total costs & \multicolumn{1}{r}{Deprivation} & \multicolumn{1}{r}{UAV} & \multicolumn{1}{r}{Truck} & \multicolumn{1}{r}{Max. deprivation} & \multicolumn{1}{r}{Demand} \\
                  & \multicolumn{1}{r}{}         & \multicolumn{1}{r}{} & \multicolumn{1}{r}{(±std.dev.)} & costs       & costs & costs                   & time (h)           & coverage                   \\ \hline
PI-bound          &   82.05  & 0.00\%  & 5293 (±199)          & 2063        & 1670  & 1560     & 18     & 75\%   \\
Re-optimization   &  944.78  & 0.00\%  & 6107 (±384)          & 2147        & 2160  & 1800     & 20     & 78\%   \\
Rule-based        &    0.01  &         & 7492 (±819)          & 1772        & 380   & 5340     & 12     & 95\%   \\
PPO               & 3535.05  &         & 6408 (±436)          & 2028        & 2220  & 2160     & 23     & 88\%   \\
DL-VFA            &  265.05  &         & 5620 (±151)          & 2100        & 1840  & 1680     & 19     & 81\%   \\
NN-VFA            & 4885.34  &         & \textbf{5585 (±213)} & 2225        & 2160  & 1200     & 20     & 75\%   \\
\hline
\multirow{2}{*}{\textbf{3 districts}} & \multicolumn{1}{r}{\multirow{2}{*}{Runtime (s)}} & \multicolumn{1}{r}{\multirow{2}{*}{Gap}} & Total costs & \multicolumn{1}{r}{Deprivation} & \multicolumn{1}{r}{UAV} & \multicolumn{1}{r}{Truck} & \multicolumn{1}{r}{Max. deprivation} & \multicolumn{1}{r}{Demand} \\
                  & \multicolumn{1}{r}{}         & \multicolumn{1}{r}{} & \multicolumn{1}{r}{(±std.dev.)} & costs       & costs & costs                   & time (h)           & coverage                   \\ \hline
PI-bound          & 1580.71  & 0.00\% & 7525 (±243)          & 2730      & 2095  & 2700     & 22     & 79\%   \\
Re-optimization   & 4713.82  & 0.01\% & 8707 (±320)          & 3287      & 3020  & 2400     & 22     & 79\%   \\
Rule-based        & 0.01     &        & 10,822 (±956)        & 2757      & 685   & 7380     & 15     & 96\%   \\
PPO               & 3299.36  &        & 10,233 (±473)        & 3328      & 3295  & 3610     & 26     & 98\%   \\
DL-VFA            & 404.46   &        & 8350 (±323)          & 3005      & 2795  & 2550     & 21     & 79\%   \\
NN-VFA            & 6768.44  &        & \textbf{7973 (±332)} & 2753      & 3150  & 2070     & 21     & 81\%   \\
\hline
\multirow{2}{*}{\textbf{4 districts}} & \multicolumn{1}{r}{\multirow{2}{*}{Runtime (s)}} & \multicolumn{1}{r}{\multirow{2}{*}{Gap}} & Total costs & \multicolumn{1}{r}{Deprivation} & \multicolumn{1}{r}{UAV} & \multicolumn{1}{r}{Truck} & \multicolumn{1}{r}{Max. deprivation} & \multicolumn{1}{r}{Demand} \\
                  & \multicolumn{1}{r}{}         & \multicolumn{1}{r}{} & \multicolumn{1}{r}{(±std.dev.)} & costs       & costs & costs                   & time (h)           & coverage                   \\ \hline
PI-bound          & 13,569.09 & 6.02\% & 10,191 (±357)          & 3111     & 2220  & 4860     & 20     & 80\%   \\
Re-optimization   & 10,780.63 & 1.14\% & 12,112 (±522)          & 4122     & 3270  & 4320     & 25     & 80\%   \\
Rule-based        & 0.01      &        & 14,300 (±1810)         & 3180     & 1280  & 9840     & 16     & 99\%   \\
PPO               & 3350.69   &        & 13,935 (±807)          & 4695     & 4710  & 4530     & 27     & 100\%  \\
DL-VFA            & 589.45    &        & 11,546 (±347)          & 4611     & 4025  & 2910     & 24     & 79\%   \\
NN-VFA            & 8796.20   &        & \textbf{11,066 (±319)} & 4156     & 4510  & 2400     & 19     & 80\%   \\
\hline
\multirow{2}{*}{\textbf{5 districts}} & \multicolumn{1}{r}{\multirow{2}{*}{Runtime (s)}} & \multicolumn{1}{r}{\multirow{2}{*}{Gap}} & Total costs & \multicolumn{1}{r}{Deprivation} & \multicolumn{1}{r}{UAV} & \multicolumn{1}{r}{Truck} & \multicolumn{1}{r}{Max. deprivation} & \multicolumn{1}{r}{Demand} \\
                  & \multicolumn{1}{r}{}         & \multicolumn{1}{r}{} & \multicolumn{1}{r}{(±std.dev.)} & costs       & costs & costs                   & time (h)            & coverage                   \\ \hline
PI-bound          & 14,404.06 & 21.53\% & 13,357 (±488)          & 4352     & 2195  & 6810     & 20     & 80\%  \\
Re-optimization   & 15,971.49 & 1.81\%  & 15,326 (±746)          & 4966     & 4780  & 5580     & 23     & 81\%  \\
Rule-based        & 0.01      &         & 16,947 (±1468)         & 4357     & 2090  & 10,500   & 18     & 98\%  \\
PPO               & 3671.79   &         & 17,542 (±1014)         & 5492     & 6500  & 5550     & 23     & 96\%  \\
DL-VFA            & 770.24    &         & 14,758 (±775)          & 4468     & 5700  & 4590     & 20     & 86\%  \\
NN-VFA            & 13,913.97 &         & \textbf{14,206 (±640)} & 5111     & 5585  & 3510     & 23     & 78\%  \\
\hline
\multirow{2}{*}{\textbf{6 districts}} & \multicolumn{1}{r}{\multirow{2}{*}{Runtime (s)}} & \multicolumn{1}{r}{\multirow{2}{*}{Gap}} & Total costs & \multicolumn{1}{r}{Deprivation} & \multicolumn{1}{r}{UAV} & \multicolumn{1}{r}{Truck} & \multicolumn{1}{r}{Max. deprivation} & \multicolumn{1}{r}{Demand} \\
                  & \multicolumn{1}{r}{}         & \multicolumn{1}{r}{} & \multicolumn{1}{r}{(±std.dev.)} & costs       & costs & costs                   & time (h)            & coverage                   \\ \hline
PI-bound          & 14,404.85  & 25.10\% & 17,699 (±729)            & 5519   & 2970  & 9210    & 22     & 79\%   \\
Re-optimization   & 20,314.81  &  2.64\% & 20,652 (±652)            & 6552   & 6040  & 8060    & 25     & 80\%   \\
Rule-based        & 0.01      &          & 23,496 (±1018)           & 5411   & 3055  & 15,030  & 19     & 98\%   \\
PPO               & 3372.81   &          & 24,784 (±848)            & 6724   & 9360  & 8700    & 28     & 98\%   \\
DL-VFA            & 948.49    &          & \textbf{19,041 (±832)}   & 6096   & 7425  & 5520    & 22     & 85\%   \\
NN-VFA            & 14,400.14  &         & 20,203 (±640)            & 6333   & 9190  & 4680    & 28     & 79\%   \\ \hline
\end{tabular}
\end{adjustbox}
\end{table}
DL-VFA provides the best computational efficiency (except for the rule-based heuristic), requiring less than 0.002s per decision epoch with 1 district and up to 0.01s per epoch for 6 districts. NN-VFA demands more computational effort with 0.03s per epoch for 1 district and 0.20s for 6 districts. These computation times are a significant improvement compared to the 15-minute computational time limit of re-optimization, a time limit in which DL-VFA can even be retrained if necessary. The solving time of the MIP models increase significantly with more districts. The training time of PPO remains stable when increasing the number of districts from 1 to 6, performing 60,000 iterations of 30 periods in less than one hour ($\pm$0.002s per decision epoch). 

\begin{figure}[ht!]
    \subfloat[3 districts]{\includegraphics[width=0.5\linewidth]{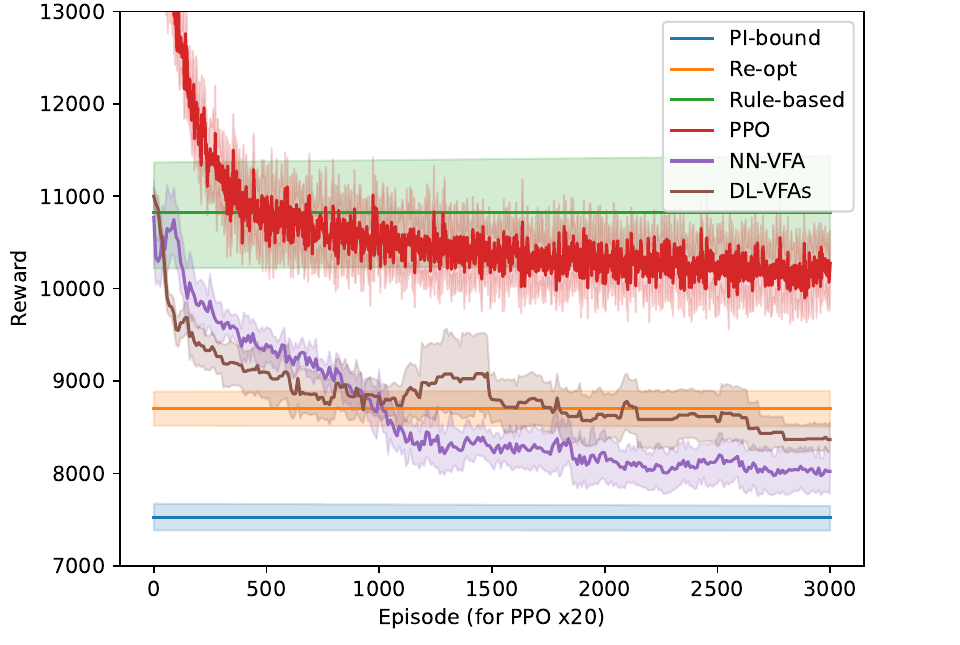}}
    \subfloat[6 districts]{\includegraphics[width=0.5\linewidth]{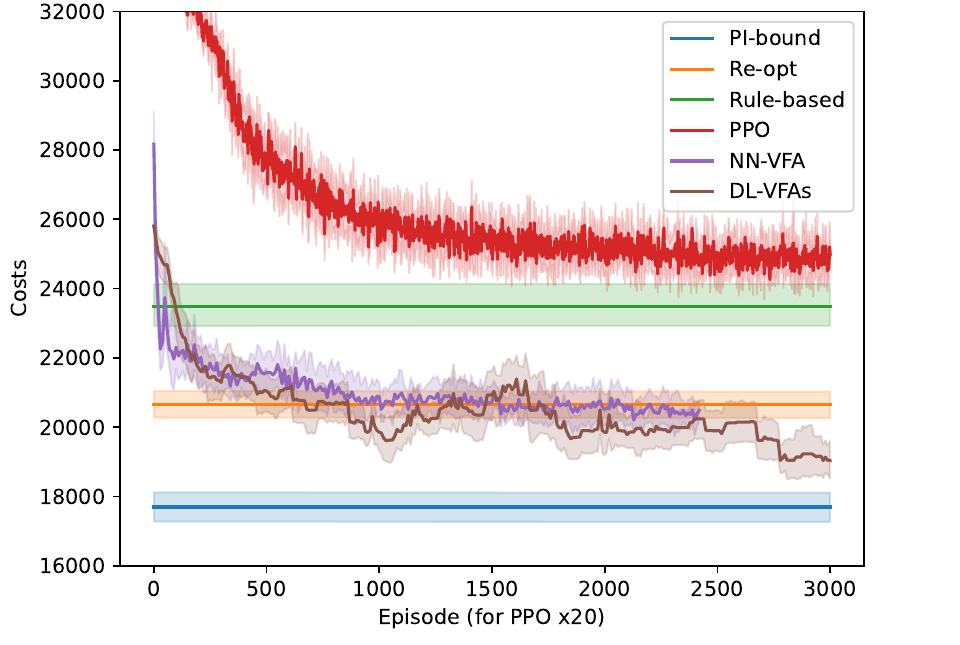}}
    \caption{Learning curves for each algorithm for 3 and 6 districts. The shaded regions around the learning curves provide the 95\% confidence intervals. Each episode is one horizon of 30 periods. For PPO, 60,000 episodes are performed and episodes in the figures are scaled to consist of 20 horizons.}
    \label{fig:learning_curves}
\end{figure}

To gain insight into the stability and training of the learning-based methods, Figure \ref{fig:learning_curves} shows the learning curves for 3 and 6 districts. The performance of DL-VFA and NN-VFA start around the costs of the rule-based heuristic and gradually decrease below the re-optimization results. For 6 districts, NN-VFA did not complete 3000 iterations within the computational time limit of 4 hours. PPO requires significantly more iterations (60,000 instead of 3000), while obtaining a performance similar to the rule-based heuristic. The learning curves for 1, 2, 4, and 5 districts can be found in Appendix \ref{app:learning_curves} and show comparable results.

\begin{figure}[ht!]
    \centering
    \subfloat[Costs for different periods and inventory levels]{\includegraphics[width=0.5\linewidth]{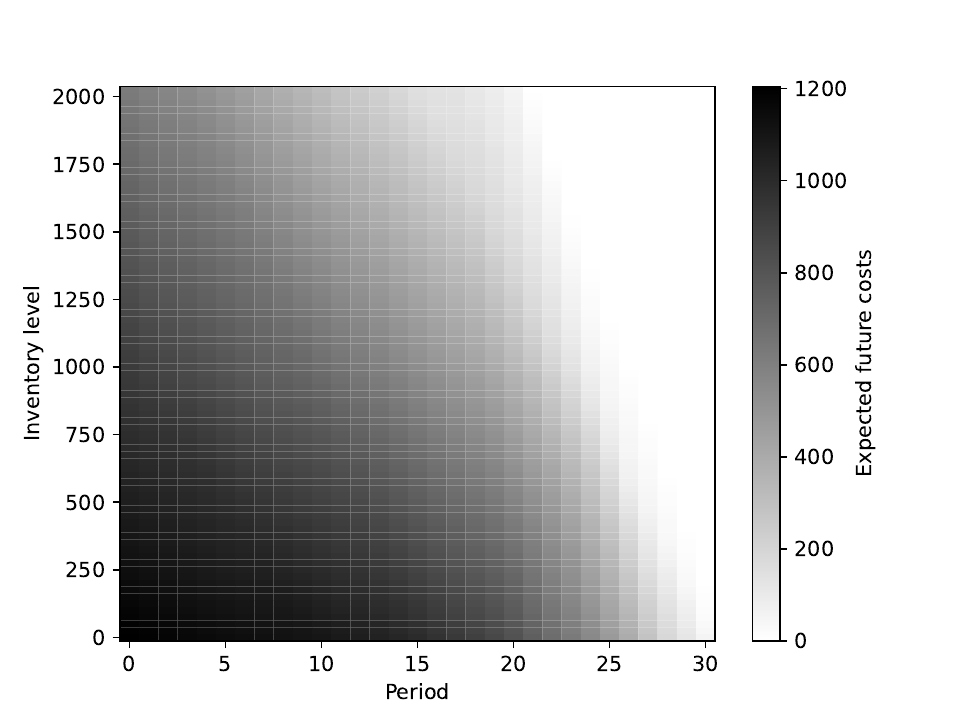}}
    \subfloat[Costs for different periods and deprivation levels]{\includegraphics[width=0.5\linewidth]{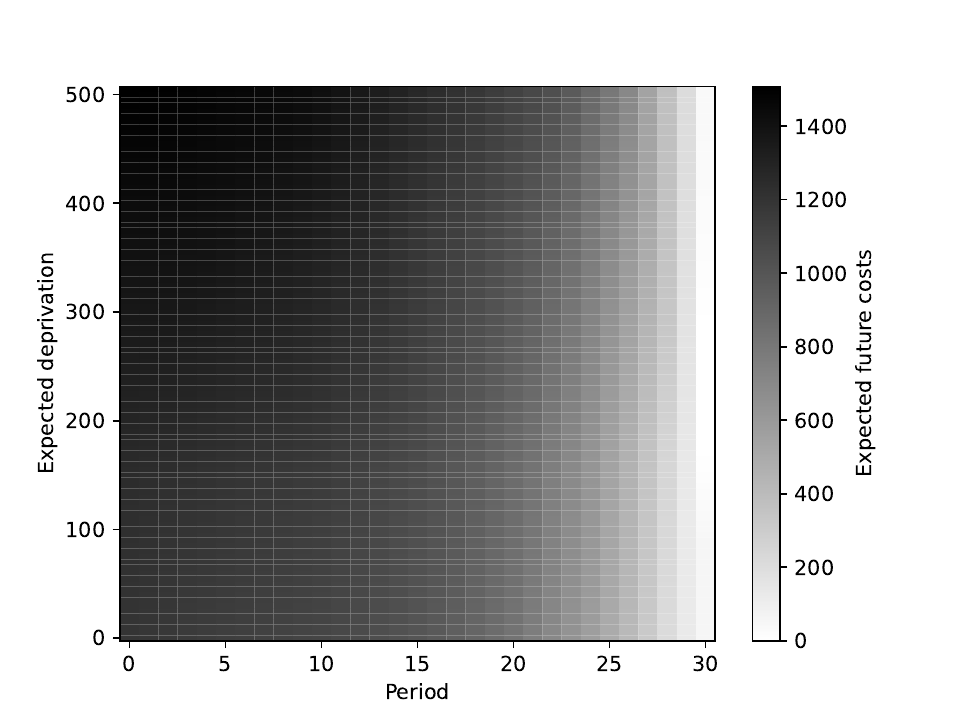}}
    \caption{Learned expected costs for different inventory levels and expected deprivation costs per period, based on the neural network trained for the 1-district instance.}
    \label{fig:dependence_plots}
\end{figure}

To provide insight into the cost effects of the neural network value function, we visualize the output of the trained neural network of NN-VFA with one district on different states in Figure \ref{fig:dependence_plots}. The future costs, predicted by the neural network, are obtained for different periods, inventory levels in the district, and expected deprivation levels. As intuitively expected, future costs decrease near the end of the horizon. Furthermore, the costs decrease significantly with higher inventory positions. Higher expected deprivation costs increase future costs, an effect that is smaller than the effect of different inventory levels.

\subsection{Results for Different Problem Variants}
\label{subsec:num_variants}
In this section, we analyze the problem with experiments featuring trucks only, different levels of uncertainty, and different patterns of demand over time. Table \ref{tab:results_variants} shows the results. DL-VFA and NN-VFA both outperform the benchmark policies in the instances with different levels of uncertainty and demand patterns, DL-VFA improves upon re-optimization by 2\% and NN-VFA improves upon re-optimization by 7.5\% on average, signifying the best results, but also the highest computational effort. The results further indicate that re-optimization performs worse under higher uncertainty, whereas DL-VFA and NN-VFA perform consistent on different levels of uncertainty. On the trucks-only instance, re-optimization performs the best. As only (larger) truckloads are allocated, the anticipation of uncertainty plays a smaller role, and optimization methods perform better than the learning-based methods. Nonetheless, comparing the trucks-only instance with the trucks and UAVs instance, the deployment of UAVs decrease deprivations costs by 22\%, transportation costs by 13\%, and maximum deprivation times from 35 to 21 hours, while the demand coverage remains 80\%. PPO performs better than the rule-based heuristic in all instance variants but does not provide competitive results.

In the scenario of decreasing demand (higher shortages initially, sufficient supplies later), both deprivation costs and UAV deployment increase. UAVs are utilized to dynamically distribute scarce supplies across the districts. In the value functions of DL-VFA, we observe larger negative weights on inventory position and larger positive weights on expected deprivation compared to other scenarios. This means that allocating goods across all districts in earlier periods significantly reduces expected future costs, which confirms the utilization of UAVs across multiple districts in case of high shortage. In contrast, in the increasing demand scenario (sufficient supplies initially, shortages later), deprivation costs are lower, and more trucks are used to allocate the larger supply at the start. In situations with sufficient supplies early on, linear value function weights on deprivation are smaller. This indicates that the focus lies on cost-efficient transport (larger shipments) while potentially neglecting another district. Spreading of supplies across all districts becomes less critical due to the higher supply rate, enabling to serve other districts the next period.
\begin{table}[ht!]
\scriptsize
\centering
\caption{Results for each method for different scenarios. The NN-VFA method mostly outperforms all other methods. Re-optimization works best for the trucks-only instance and well for the low-uncertainty variant. The gap for re-optimization is the average gap over the 30 solutions in the horizon.}
\label{tab:results_variants}
\begin{adjustbox}{center}
\begin{tabular}{lrrrrrrrr}
\hline
\multirow{2}{*}{\textbf{Trucks only}} & \multicolumn{1}{r}{\multirow{2}{*}{Runtime (s)}} & \multicolumn{1}{r}{\multirow{2}{*}{Gap}} & Total costs                     & \multicolumn{1}{r}{Deprivation} & \multicolumn{1}{r}{UAV} & \multicolumn{1}{r}{Truck} & \multicolumn{1}{r}{Max. deprivation} & \multicolumn{1}{r}{Demand} \\
                  & \multicolumn{1}{r}{}         & \multicolumn{1}{r}{} & \multicolumn{1}{r}{(±std.dev.)} & costs       & costs & costs                   & time (h)           & coverage                   \\ \hline
PI-bound        & 1113.18 & 0.00\% & 8788 (±333)          & 3358 & 0 & 5430 & 35 & 80\% \\
Re-optimization & 2968.10 & 0.01\% & \textbf{9500 (±573)} & 3510 & 0 & 5990 & 35 & 80\% \\
Rule-based      & 0.00    &        & 12,191 (±735)        & 3251 & 0 & 8940 & 24 & 97\% \\
PPO             & 3647.49 &        & 10,893 (±826)        & 3243 & 0 & 7650 & 30 & 99\% \\
DL-VFA          & 342.30  &        & 10,045 (±614)        & 3855 & 0 & 6190 & 37 & 85\% \\
NN-VFA          & 7145.99 &        & 9923 (±250)          & 3223 & 0 & 6700 & 36 & 87\% \\
\hline
\multirow{2}{*}{\textbf{Low uncertainty}} & \multicolumn{1}{r}{\multirow{2}{*}{Runtime (s)}} & \multicolumn{1}{r}{\multirow{2}{*}{Gap}} & Total costs & \multicolumn{1}{r}{Deprivation} & \multicolumn{1}{r}{UAV} & \multicolumn{1}{r}{Truck} & \multicolumn{1}{r}{Max. deprivation} & \multicolumn{1}{r}{Demand} \\
                  & \multicolumn{1}{r}{}         & \multicolumn{1}{r}{} & \multicolumn{1}{r}{(±std.dev.)} & costs       & costs & costs                   & time (h)            & coverage                   \\ \hline
PI-bound        & 1410.15 & 0.00\% & 7570 (±147)          & 2820 & 2200 & 2550 & 20 & 78\%   \\
Re-optimization & 3767.23 & 0.05\% & 8554 (±118)          & 3459 & 2755 & 2340 & 25 & 77\%   \\
Rule-based      & 0.00    &        & 11,357 (±986)        & 2712 & 665  & 7980 & 16 & 96\%   \\
PPO             & 3697.48 &        & 10,041 (±315)        & 2881 & 3520 & 3640 & 25 & 96\%   \\
DL-VFA          & 417.39  &        & 8506 (±564)          & 2931 & 3145 & 2430 & 24 & 81\%   \\
NN-VFA          & 7596.05 &        & \textbf{8016 (±126)} & 2706 & 3480 & 1830 & 21 & 81\%   \\
\hline
\multirow{2}{*}{\textbf{High uncertainty}} & \multicolumn{1}{r}{\multirow{2}{*}{Runtime (s)}} & \multicolumn{1}{r}{\multirow{2}{*}{Gap}} & Total costs & \multicolumn{1}{r}{Deprivation} & \multicolumn{1}{r}{UAV} & \multicolumn{1}{r}{Truck} & \multicolumn{1}{r}{Max. deprivation} & \multicolumn{1}{r}{Demand} \\
                  & \multicolumn{1}{r}{}         & \multicolumn{1}{r}{} & \multicolumn{1}{r}{(±std.dev.)} & costs       & costs & costs                   & time (h)            & coverage                   \\ \hline
PI-bound        & 2015.60 & 0.00\% & 7421 (±371)          & 2536 & 1765 & 3120 & 19 & 80\% \\
Re-optimization & 5610.57 & 0.12\% & 8921 (±497)          & 2951 & 3420 & 2550 & 23 & 82\% \\
Rule-based      & 0.01    &        & 10,824 (±1189)       & 2624 & 640  & 7560 & 16 & 97\% \\
PPO             & 3307.25 &        & 10,298 (±808)        & 2998 & 3570 & 3730 & 26 & 95\% \\
DL-VFA          & 387.43  &        & 8542 (±284)          & 2382 & 3490 & 2670 & 26 & 87\% \\
NN-VFA          & 6463.33 &        & \textbf{8036 (±384)} & 2746 & 3190 & 2100 & 22 & 82\% \\
\hline
\multirow{2}{*}{\textbf{Decreasing demand}} & \multicolumn{1}{r}{\multirow{2}{*}{Runtime (s)}} & \multicolumn{1}{r}{\multirow{2}{*}{Gap}} & Total costs & \multicolumn{1}{r}{Deprivation} & \multicolumn{1}{r}{UAV} & \multicolumn{1}{r}{Truck} & \multicolumn{1}{r}{Max. deprivation} & \multicolumn{1}{r}{Demand} \\
                  & \multicolumn{1}{r}{}         & \multicolumn{1}{r}{} & \multicolumn{1}{r}{(±std.dev.)} & costs       & costs & costs                   & time (h)            & coverage                   \\ \hline
PI-bound            & 2166.06 & 0.00\% & 8624 (±351)          & 3854 & 2460 & 2310 & 19 & 67\% \\
Re-optimization     & 4876.01 & 0.27\% & 10,048 (±407)        & 3933 & 3715 & 2400 & 25 & 71\% \\
Rule-based          & 0.00    &        & 11,880 (±785)        & 3475 & 1055 & 7350 & 13 & 91\% \\
PPO                 & 3314.15 &        & 11,389 (±501)        & 4114 & 4375 & 2900 & 28 & 92\% \\
DL-VFA              & 515.10  &        & 9811 (±330)          & 3321 & 4990 & 1500 & 23 & 74\% \\
NN-VFA              & 7049.12 &        & \textbf{9547 (±422)} & 3572 & 4355 & 1620 & 21 & 71\% \\
\hline
\multirow{2}{*}{\textbf{Increasing demand}} & \multicolumn{1}{r}{\multirow{2}{*}{Runtime (s)}} & \multicolumn{1}{r}{\multirow{2}{*}{Gap}} & Total costs & \multicolumn{1}{r}{Deprivation} & \multicolumn{1}{r}{UAV} & \multicolumn{1}{r}{Truck} & \multicolumn{1}{r}{Max. deprivation} & \multicolumn{1}{r}{Demand} \\
                  & \multicolumn{1}{r}{}         & \multicolumn{1}{r}{} & \multicolumn{1}{r}{(±std.dev.)} & costs       & costs & costs                   & time (h)           & coverage                   \\ \hline
PI-bound        & 2596.36 & 0.00\% & 6347 (±249)          & 1602 & 845  & 3900 & 22 & 90\% \\
Re-optimization & 5160.64 & 0.12\% & 7505 (±331)          & 1775 & 1980 & 3750 & 25 & 92\% \\
Rule-based      & 0.00    &        & 9373 (±856)          & 2053 & 510  & 6810 & 12 & 96\% \\
PPO             & 3695.31 &        & 8894 (±671)          & 2804 & 2860 & 3230 & 21 & 98\% \\
DL-VFA          & 430.07  &        & 7449 (±303)          & 2774 & 2365 & 2310 & 30 & 86\% \\
NN-VFA          & 6193.41 &        & \textbf{6828 (±149)} & 1888 & 2660 & 2280 & 20 & 93\% \\
\hline
\end{tabular}
\end{adjustbox}
\end{table}

\begin{figure}[ht!]
    \centering
     \makebox[\textwidth][c]{
        \subfloat[Trucks only (Deprivation: 3484 - Transport: 6900)]{\includegraphics[width=0.63\linewidth]{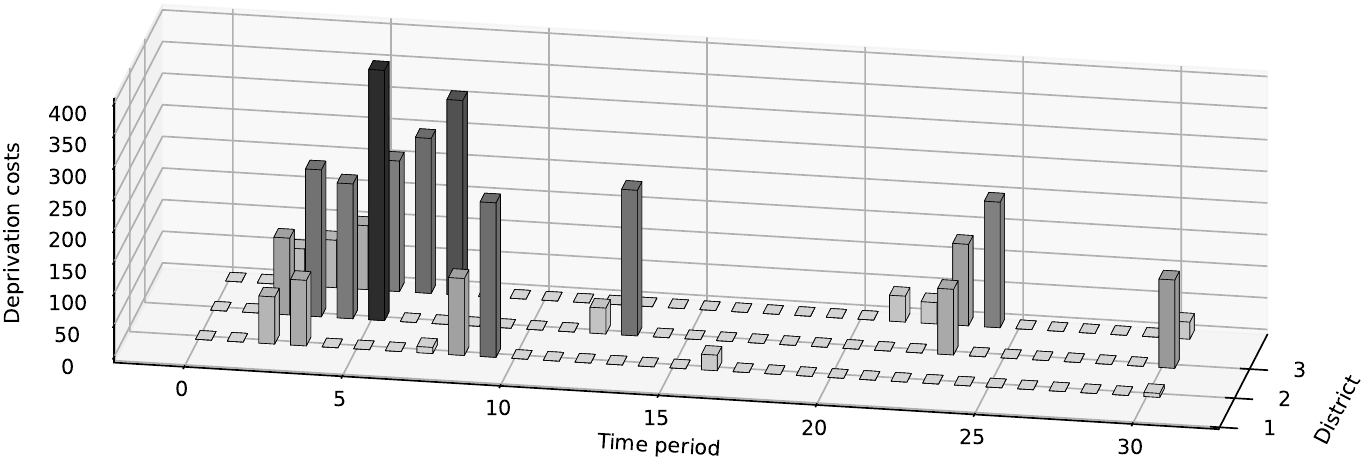}}
        \subfloat[Trucks \& UAVs (Deprivation: 2528 - Transport: 5300)]{\includegraphics[width=0.63\linewidth]{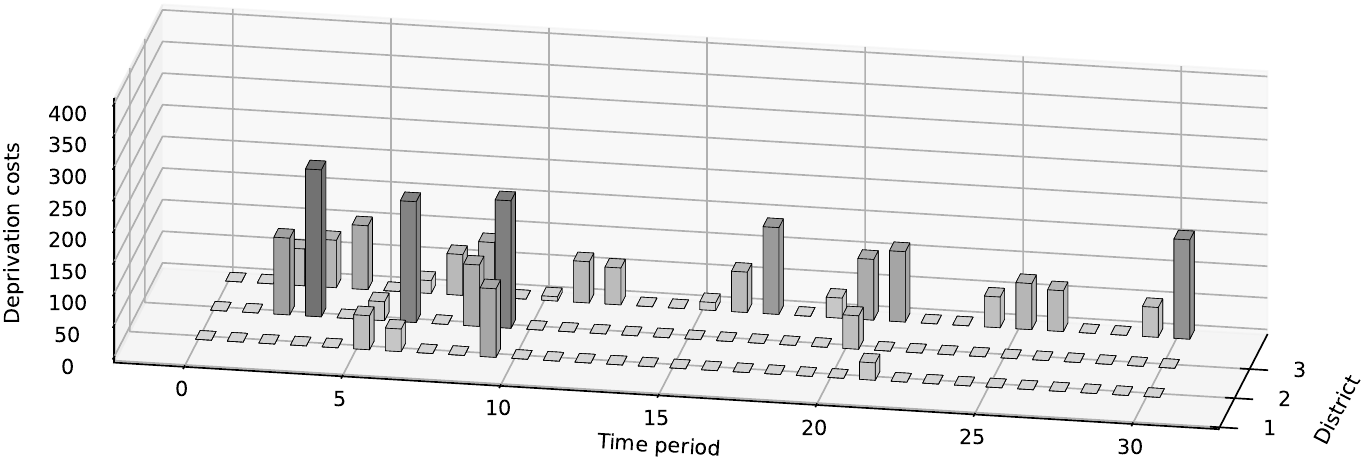}}
    }
    \caption{Deprivation costs over time per district, for one realization of the 3-district instance. In the trucks-only instance, deprivation is high at the start in multiple districts, whereas in the trucks and UAVs instance, transportation- and deprivation costs are overall lower and more spread over districts and time.}
    \label{fig:deprivation_time_td}
\end{figure}

To gain more insight into the added value of UAVs compared to a trucks-only scenario, a spatio-temporal overview of deprivation costs is visualized in Figure \ref{fig:deprivation_time_td} for three districts. In the trucks-only scenario, deprivation costs mainly arise in the early stages. With trucks and UAVs, the deprivation is overall lower, especially in the early stages, and also more spread over time and districts. Note that the number of consecutive periods of deprivation is lower for the trucks-and-UAV scenario than for the trucks-only scenario. Districts with relatively high demand (districts 1 and 2 in Figure \ref{fig:deprivation_time_td}) receive goods by truck, whereas districts with lower demand (district 3) receive supplies by UAV. UAV deliveries can be identified by the regular reset in deprivation in district 3. Besides a reduction in deprivation costs of 22\% and a reduction of transportation costs of 13\%, the deployment of UAVs also provides a more equal distribution of the deprivation costs over time and districts. In this way, UAVs aid in a better trade-off between the egalitarian perspective (equal service to everyone) and the utilitarian perspective (focus on the largest groups to deliver the most).

\subsection{Nepal Case}
\label{subsec:num_Nepal}
With the Nepal case, we analyze the suitability of the proposed methods on a practice-based instance, both in performance and computational time required to solve the larger-scale problem. The experiment is based on the Nepal 2015 Earthquake, where 13 districts and the capital Kathmandu were heavily damaged. We analyze the same metrics as with the smaller instances, and compare the scenario without UAVs (as occurred in 2015) and with UAVs to estimate the potential benefit of UAVs in practice. Furthermore, we analyze the decisions over time in terms of truck and UAV allocations to study their added value in different stages of the operation.

The computational results of the Nepal case are presented in Table \ref{tab:results_Nepal} for both the trucks-only scenario and the scenario with trucks and UAVs. The optimization gaps in the MIP models increase compared to the smaller instances of Section \ref{subsec:num_districts}. With a training time of less than an hour, DL-VFA provides the best results, improving upon re-optimization by 8\%, while NN-VFA improves upon re-optimization by 5\%. DL-VFA and NN-VFA improve the rule-based heuristic and PPO even more. These results show that although PPO is suitable for continuous action spaces, it does not provide fruitful results for this combinatorial action space. Similar to the smaller trucks-only scenario, re-optimization finds the best results in the trucks-only instance of Nepal. The optimization methods are likely to perform better than the learning-based methods as only larger truckloads are allocated, covering demand over multiple periods and thereby mitigating the effects of uncertainty. Nonetheless, the trucks and UAV instance, compared to the trucks-only instance, reduces both deprivation costs and transportation costs by 18\% and reduces the maximum deprivation time across all districts from 36 hours to 26 hours.
\begin{table}[ht!]
\scriptsize
\centering
\caption{Results for each method in the Nepal case. The gap for re-optimization is the average gap over the 30 solutions in the horizon.}
\label{tab:results_Nepal}
\begin{adjustbox}{center}
\begin{tabular}{lrrrrrrrr}
\hline
\multirow{2}{*}{\textbf{Trucks \& UAVs}} & \multicolumn{1}{r}{\multirow{2}{*}{Runtime (s)}} & \multicolumn{1}{r}{\multirow{2}{*}{Gap}} & Total costs                     & \multicolumn{1}{r}{Deprivation} & \multicolumn{1}{r}{UAV} & \multicolumn{1}{r}{Truck} & \multicolumn{1}{r}{Max. deprivation} & \multicolumn{1}{r}{Demand} \\
                  & \multicolumn{1}{r}{}         & \multicolumn{1}{r}{} & \multicolumn{1}{r}{(±std.dev.)} & costs       & costs & costs                   & time (h)            & coverage                   \\
\hline
PI-bound              & 14,403.21 & 38.51\% & 29,660 (±727)           & 8941   & 6852   & 13,867 & 21 & 86\% \\
Re-optimization       & 23,588.47 & 4.49\%  & 33,130 (±1230)          & 9439   & 9844   & 13,847 & 25 & 89\% \\
Rule-based            & 0.01      &         & 34,200 (±2016)          & 12,555 & 8179   & 13,466 & 18 & 97\% \\
PPO                   & 4023.88   &         & 44,303 (±716)           & 10,016 & 34,286 & 0      & 24 & 99\% \\
DL-VFA                & 3282.22   &         & \textbf{30,464 (±1031)} & 8816   & 12,298 & 9350   & 26 & 91\% \\
NN-VFA                & 14,404.18 &         & 31,336 (±1047)          & 10767  & 11,228 & 9342   & 32 & 90\% \\              
\hline
\multirow{2}{*}{\textbf{Trucks only}} & \multicolumn{1}{r}{\multirow{2}{*}{Runtime (s)}} & \multicolumn{1}{r}{\multirow{2}{*}{Gap}} & Total costs & \multicolumn{1}{r}{Deprivation} & \multicolumn{1}{r}{UAV} & \multicolumn{1}{r}{Truck} & \multicolumn{1}{r}{Max. deprivation} & \multicolumn{1}{r}{Demand} \\
                  & \multicolumn{1}{r}{}         & \multicolumn{1}{r}{} & \multicolumn{1}{r}{(±std.dev.)} & costs       & costs & costs                   & time (h)           & coverage                    \\
\hline
PI-bound              & 14,404.61 & 32.92\% & 35,819 (±1201)          & 9880   & 0   & 25,939  & 31   & 86\% \\
Re-optimization       & 20,819.06 & 2.61\%  & \textbf{37,172 (±1328)} & 10,744 & 0   & 26,427  & 32   & 86\% \\
Rule-based            & 0.01      &         & 92,423 (±2985)          & 73,328 & 0   & 19,095  & 84   & 98\% \\
PPO                   & 4216.09   &         & 51,958 (±3069)          & 24,351 & 0   & 27,067  & 52   & 99\% \\
DL-VFA                & 2328.06   &         & 40,094 (±2068)          & 15,881 & 0   & 24,213  & 48   & 94\% \\
NN-VFA                & 14,403.11 &         & 50,589 (±4077)          & 13,958 & 0   & 36,631  & 42   & 99\% \\
\hline
\end{tabular}
\end{adjustbox}
\end{table}
When analyzing the results at the district level, we find that districts with relatively lower transportation costs and higher demands receive relatively more supplies, resulting in lower deprivation costs and higher demand coverage than districts with lower demands and higher transportation costs. In this sense, the consideration of deprivation costs aids in balancing the trade-off between the egalitarian perspective and the utilitarian perspective. We also found, similar to Figure \ref{fig:deprivation_time_td}, that the use of UAVs results in a more even spatio-temporal spread of allocations and deprivations, as reflected by the lower maximum deprivation times presented in Table \ref{tab:results_Nepal}.
 
\begin{figure}[ht!]
    \centering
    \includegraphics[width=0.43\linewidth]{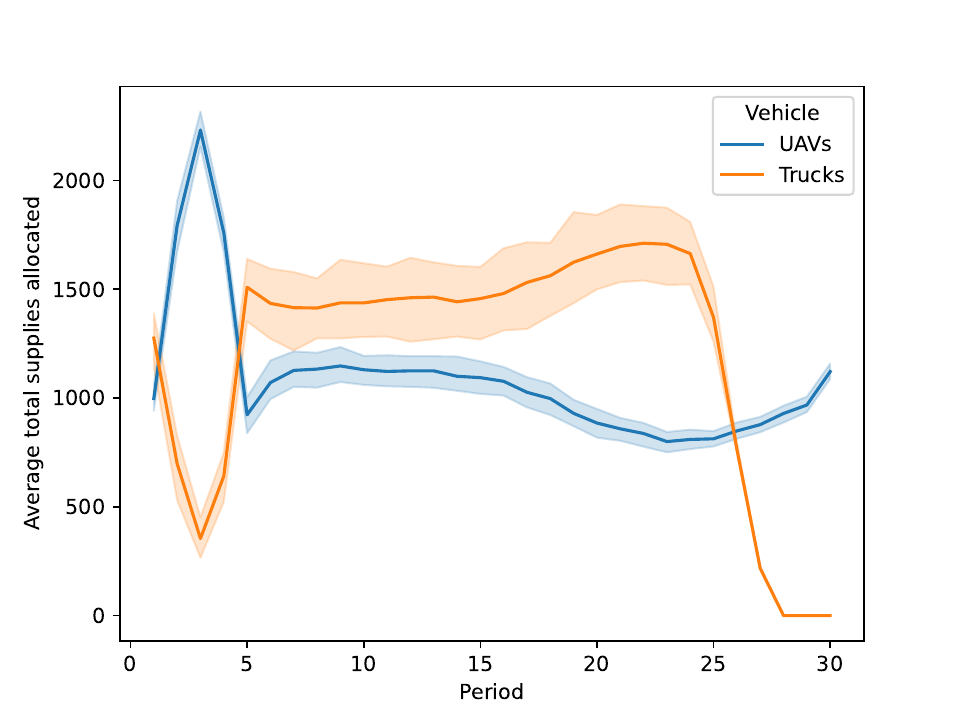}
    \caption{Sum of allocated supplies per period per vehicle in the Nepal case with trucks and UAVs.} \label{fig:allocations_per_vehicle_Nepal}
\end{figure}

To analyze truck- and UAV allocations over time, we averaged the number of supplies allocated per vehicle per decision epoch over 1000 episodes. Figure \ref{fig:allocations_per_vehicle_Nepal} shows the pattern that emerges over time. Note that we assumed a stable demand pattern (i.e., no increase or decrease over time). Especially the beginning and end of the horizon show interesting effects. In the first (and second) period, supplies are usually allocated using trucks to a few districts with high demand. Then, generally in the second, third, and fourth period, supplies are allocated to all remaining districts using UAVs to alleviate deprivation everywhere. The situation stabilizes thereafter and the majority of supplies are delivered by trucks, and mostly low-demand districts are visited by UAVs. Truck deliveries slightly increase by more available supplies over time. This indicates that for the long-term and large quantities, UAVs are less suitable. Near the end, truck allocations decline (as the remaining demand is limited) and remaining deliveries are completed with UAVs.

To summarize all experiments, in the instances with varying numbers of districts, we find that DL-VFA provides the best computational efficiency and competitive results. The Nepal case also reflects these results, in which DL-VFA improves the benchmarks by 8\%-31\% with limited training time. Moreover, for the smaller instances (up to 5 districts), NN-VFA would be the recommended approach and NN-VFA also adapts well to different problem variants. With abundant computational resources, this method may also be a compelling option for larger instances. For larger instances, better scalability, and explainability through linear weights, DL-VFA would be the preferred method providing excellent performance. We also highlight the significant role of UAVs, not only found in different problem variants but also in the Nepal case, in which UAVs are essential in the initial and final stages to distribute scarce supplies evenly or deliver smaller amounts relatively cheaply, reducing deprivation costs and transportation costs by 18\% compared to using only trucks. Lastly, the consideration of deprivation costs in combination with UAVs results in a reasonably even spatio-temporal distribution of supplies and deprivation, providing an efficient operation while striking a balance between egalitarian and utilitarian relief.

\section{Conclusions}
\label{sec:conclusions}
In this paper, we investigated the impact of deprivation costs and the role of UAVs by introducing the stochastic dynamic post-disaster inventory allocation problem (SDPDIAP) with uncertain supply and demand. Deprivation costs are an important metric to address the spatio-temporal consequences of allocating scarce relief supplies over multiple affected districts, whereas humanitarian UAVs offer the potential to cost-effectively and flexibly distribute supplies in the aftermath of a disaster. We address deprivation costs in a stochastic dynamic setting and formulate the SDPDIAP as a Markov decision process model. As this model is computationally intractable to solve, we propose two anticipatory policies based on approximate dynamic programming. First, we propose a decomposed linear value function approximation (DL-VFA) that decomposes the problem in spatio-temporal value functions, providing high performance and excellent scalability, with the combinatorial allocation of supplies in each period being handled by solving an mixed integer program (MIP) that includes the spatially decomposed value functions of the respective period. Second, we propose a neural network value function approximation (NN-VFA) that captures the spatio-temporal and nonlinear effects in a single neural network, which is integrated into an MIP that is used to solve the single-stage combinatorial allocation problem. We compare the performance of DL-VFA and NN-VFA with a perfect information bound, an exact re-optimization approach, a rule-based heuristic, and proximal policy optimization (PPO). Up to 5 disticts, NN-VFA provides the best results (2.6\% less costs than DL-VFA and 8.3\% less costs than re-optimization on average). DL-VFA shows excellent scalability, with the ability to train value functions with limited computational effort, improving upon NN-VFA by 4.3\% and re-optimization by 7\% for larger instances (6 and 13 districts). Furthermore, DL-VFA provides clear explainability by the linear weights of the value functions, which is a valuable trait in disaster environments. DL-VFA improves the rule-based heuristic and PPO on average by 22.7\% and 23.7\%, respectively, whereas NN-VFA improves the rule-based heuristic and PPO on average by 16.9\% and 17.5\%, respectively. NN-VFA proves to be the preferred approach for smaller instances and in situations with ample computational resources, whereas DL-VFA is preferred for larger instances, situations requiring computational efficiency and explainability of decisions.

By considering both transportation- and deprivation costs, the best solutions allocate goods to all districts, but districts with higher demands and lower transportation costs receive relatively more. This indicates a balance between the egalitarian perspective (equal service to everyone) and the utilitarian perspective (focus on the largest groups). Results show that districts with relatively higher transportation costs and lower demands are visited more by UAVs than trucks. When comparing trucks-only scenarios with scenarios including trucks and UAVs, the addition of UAVs provides reductions in deprivation costs of 18-22\%, reductions in transportation costs between 13-18\%, and reductions in maximum deprivation times of 19-40\%, while maintaining a similar coverage of demand. Generally speaking, multiple transportation modes, e.g., a bulk mode and a flexible mode (in our case trucks and UAVs) broaden the opportunities for humanitarian operations that are both cost-efficient and effective in alleviating human suffering. UAVs especially play a crucial role at the operations' start and end by distributing scarce supplies in smaller amounts among various districts.

There are several opportunities for future work. The spatio-temporal decomposition allows for extensions in decision-making, such as (i) a multi-depot problem, in which supplies from multiple warehouses are allocated to districts, or (ii) inter-district transhipments, with supplies being rebalanced between districts after initial allocation. It might be interesting to address accessibility constraints due to disruptions in the infrastructure or to extend the model towards multiple types of relief supplies. Furthermore, as the single-stage decision problems are solved with MIPs, future endeavors may add practical constraints or operating rules for decision-making.

\bibliographystyle{informs2014trsc}
\bibliography{references}

\clearpage

\appendix

\section{Literature Overview}
\label{app:literature}

\begin{table}[ht!]
\scriptsize
\caption{Overview of the literature related to our study}
\label{tab:literature}
\begin{adjustbox}{center}
\begin{tabular}{p{3.5cm}p{2.5cm}p{2.1cm}p{2.6cm}p{2.0cm}p{1.9cm}p{1.7cm}}
\hline
Reference                                          & Problem                       & Sources of uncertainty        & Objective                          & Modality                   & Decomposition          & Approach               \\
\hline
\cite{ben2011robust}              & Evacuee flow assignment       & Evacuee demand                & Travel + penalty costs             & -                          & -                      & ARO                    \\
\cite{vanajakumari2016integrated} & Location-allocation-routing   & Demand                        & Response times                     & Three truck types          & -                      & MIP                    \\
\cite{perez2016inventory}         & Allocation-routing            & -                             & Transport + deprivation costs      & Trucks                     & -                      & MIP                    \\
\cite{rivera2016dynamic}          & Relief allocation             & -                             & Transport + deprivation costs      & Trucks, boats              & Time                   & MPC                    \\
\cite{nadi2016reinforcement}      & Team allocation               & Demand, road damage           & Assessment time                    & Teams                      & Time                   & Q-learning             \\
\cite{bouzaiene2016single}        & Assignment                    & Travel times, delays          & Costs                              & Trains                     & Time, locomotives      & ADP                    \\
\cite{nadi2017adaptive}           & Team allocation               & Demand, road damage           & Assessment time                    & Two types of teams         & Time                   & Q-learning             \\
\cite{natarajan2017multi}         & Inventory-allocation          & Funding, demand               & Costs                              & -                          & Time                   & Rule-based heuristic   \\
\cite{rivera2017anticipatory}     & Freight selection             & Demand                        & Costs                              & Trucks and barges          & Time                   & ADP                    \\
\cite{moreno2018effective}        & Location-routing              & Demand, supply, access & Transport + deprivation costs & Trucks, boats, helicopters & -                      & Two-stage SMIP         \\
\cite{yu2019rollout}              & Relief allocation             & -                             & Transport + deprivation costs      & -                          & Time                   & Rollout algorithm      \\
\cite{liu2019robust}              & Distribution                  & Demand, supply                & Unmet demand                       & Helicopters                & Time                   & Robust MPC             \\
\cite{zhan2021disaster}           & Relief allocation             & Demand                        & Unmet demand                       & Two truck types            & Time                   & PSO                    \\
\cite{yu2021reinforcement}        & Relief allocation             & -                             & Transport + deprivation costs      & -                          & Time                   & Q-learning             \\
\cite{fan2022dhl}                 & Relief allocation             & -                             & Transport + deprivation costs      & -                          & Time                   & Deep Q-learning        \\
\cite{ghasbeh2022equitable}       & Location-distribution         & Demand, budget                & Unmet demand + travel time         & Road vehicle               & -                      & ALNS/MDLS              \\
\cite{zhang2022online}            & Parking allocation            & Demand                        & Travel time                        & -                          & Parking areas          & Multi-agent Q-learning \\
\cite{beirigo2022learning}        & Mobility on demand            & Demand                        & Profit                             & Owned and hired vehicles    & Time, area, vehicle type     & ADP                    \\ 
\textbf{Our paper}                        & \textbf{Relief allocation}             & \textbf{Demand, supply}                & \textbf{Transport + deprivation costs}    & \textbf{Trucks, UAVs}              & \textbf{Time, districts}            & \textbf{ADP} \\
\hline
\end{tabular}
\end{adjustbox}
\end{table}

\clearpage

\section{Warm-up heuristic}
\label{app:warm-up}
\begin{figure}[!ht]
\centering
\begin{minipage}{1.0\linewidth}
\begin{algorithm}[H]
    \footnotesize
    \linespread{1}\selectfont
    \caption{Warm-up heuristic}
    \label{alg:warmup}
    \begin{algorithmic}
        \State \textbf{Step 0.} Set buffer size $B$ and discount rate $\lambda$
        \For{$b = 1, 2, \dots, B$}
            \State \textbf{Step 1.} Generate sample path $\omega^b$.
            \State \textbf{Step 2.} Set the initial state to $S_0$ using $W_0 \in \omega^b$.
            \For{$t = 0, 1,..., T $}
                \State \textbf{Step 3a.} Obtain decision $\mathbf{x}_t$:
                    \State $I^{\text{CW,x}}_t = I^{\text{CW}}_t$
                    \For {$n \in \set{N}$}
                        \If{$\delta_{tn} \geqslant Z_1 \sim \textnormal{unif}\{1,3\}$}
                            \State $x^{\text{interim}}_{tn} = q_0 \cdot Z_2 \sim \text{unif}\{1,3\}$
                            \If{$x^{\text{interim}}_{tn} \geqslant I^{\text{CW,x}}_t$}
                                \State $x_{tn0} = x^{\text{interim}}_{tn}$
                                \State $I^{\text{CW,x}}_t \mathrel{-}= x^{\text{interim}}_{tn}$
                            \EndIf
                        \EndIf
                    \EndFor
                    \State $n_1 = Z_3~\text{unif}\{1,N\}$
                    \State $x^{\text{interim}}_{tn} = \min\left\{q_1, I^{\text{CW,x}}_t + x_{t,n_1,0} \right\}$
                    \If{$\delta_{t, n_1} \geqslant Z_3 \sim \textnormal{unif}\{1,3\}$}
                        \State  $x_{tn0} = 0$
                        \State $x_{tn1} = x^{\text{interim}}_{tn}$
                    \EndIf
                \StateLong{\textbf{Step 3b.} Compute and store post-decision states $S^\text{x}_t = S^{Mx}(S_t, \mathbf{x}_t)$ and related costs $C(S_t, \mathbf{x}_t)$ in buffer $\set{B}^b$.}
                \StateLong{\textbf{Step 3c.} Obtain new information $W_{t+1} \in \omega^b$ and transition to the next state $S_{t+1} = S^{MW}(S^\text{x}_t, W_{t+1})$.}
            \EndFor
            \StateLong{\textbf{Step 4a.} Remove outliers from buffer $\set{B}^r$ if direct costs in sample path are higher than $Q3 + 1.5IQR$.}
            \StateLong{\textbf{Step 4b.} Recursively compute realized future district costs $\hat{V}_{tn}$ with algorithmic discount factor $\lambda \quad \forall \enspace t \in \set{T}, \enspace n \in \set{N}$.}
        \EndFor
        \State \textbf{Step 5.} \textbf{return} experience buffer $\set{B}$ 
    \end{algorithmic}
\end{algorithm}
\end{minipage}
\end{figure}

\clearpage

\section{Neural network value function approximation}
\label{app:NN-VFA}
\begin{figure}[!ht]
\centering
\begin{minipage}{1.0\linewidth}
\begin{algorithm}[H]
    \footnotesize
    \linespread{1}\selectfont
    \caption{Neural network VFA algorithm}
    \label{alg:NN-VFA}
    \begin{algorithmic}
        \State \textbf{Step 0.} Set buffer size $B$ and discount factor $\lambda$.
        \StateLong{\textbf{Step 1.} Initialize exploration parameter $\epsilon^0$, update frequency $u$, and buffer $\set{B}^0$ for $B$ horizons to fill the experience buffer by running the warm-up heuristic.}
        \StateLong{\textbf{Step 2.} Compute neural network $f(\Phi^{\text{neural}}(S^\text{x}_t) | \Theta^{\text{neural},r})$ based on $\set{B}^0$.}
        \For{$r = 1, 2, \dots, R$}
            \State \textbf{Step 3.} Generate sample path $\omega^r$.
            \State \textbf{Step 4.} Set the initial state to $S_0$ using $W_0 \in \omega^r$.
            \For{$t = 0, 1,..., T $}
                \State \textbf{Step 5a.} Obtain decision $\mathbf{x}_t$:
                    \If{$Z \sim \textnormal{unif}(0,1) < \epsilon^r$}
                        \State Get $\mathbf{x}_t$ by using step 3a from the warm-up heuristic.
                    \Else
                        \StateLong{Get $\mathbf{x}_t$ by solving the optimization problem (\ref{eq:nn_obj})-(\ref{eq:nn_misc}) using neural network $\mathbb{F}(\Phi^{\text{neural}}(S^\text{x}_t) | \Theta^{\text{neural}})$.}
                    \EndIf
                \StateLong{\textbf{Step 5b.} Compute and store post-decision states $S^\text{x}_t = S^{Mx}(S_t, \mathbf{x}_t)$ and related costs $C(S_t, \mathbf{x}_t)$ in buffer $\set{B}^r$.}
                \StateLong{\textbf{Step 5c.} Obtain new information $W_{t+1} \in \omega^r$ and transition to the next state $S_{t+1} = S^{MW}(S^\text{x}_t, W_{t+1})$.}
            \EndFor
            \If{$r$ is multiple of $u$}
                \StateLong{\textbf{Step 6a.} Remove outliers from buffer $\set{B}^r$ if direct costs in sample path are higher than $Q3 + 1.5IQR$.}
                \StateLong{\textbf{Step 6b.} Recursively compute realized future district costs $\hat{V}_{t}$ with algorithmic discount factor $\lambda \quad \forall \enspace t \in \set{T}$.}
                \StateLong{\textbf{Step 6c.} Update weights $\hat{\Theta}^{\text{neural},r}$ of the neural network value function by backpropagation through the network.}
                \State \textbf{Step 6d.} Update $\epsilon^r$.
            \EndIf
        \EndFor
        \State \textbf{Step 7. return} value function $\overline{V}^R$.
    \end{algorithmic}
\end{algorithm}
\end{minipage}
\end{figure}

\clearpage

\section{Instance Data}
\label{app:instances}
Details on the instances with different districts are presented in Table \ref{tab:app:details_districts}, the specific demand patterns of the decreasing and increasing demand scenarios in Table \ref{tab:app:demand_patterns_variants}, and the details of the Nepal case are shown in Table \ref{tab:app:details_Nepal_instance}.
\begin{table}[ht!]
    \centering
    \footnotesize
    \caption{Demand per district per period for the increasing demand scenario and decreasing demand scenario}
    \label{tab:app:details_districts}
\begin{tabular}{llll}
\hline
District & Demand per period & UAV costs & Truck costs \\ \hline
1        & 200               & 150       & 900         \\ \hline
District & Demand per period & UAV costs & Truck costs \\ \hline
1        & 300               & 100       & 600         \\
2        & 100               & 200       & 1200        \\ \hline
District & Demand per period & UAV costs & Truck costs \\ \hline
1        & 200               & 50        & 300         \\
2        & 300               & 150       & 900         \\
3        & 100               & 250       & 1500        \\ \hline
District & Demand per period & UAV costs & Truck costs \\ \hline
1        & 200               & 50        & 300         \\
2        & 300               & 150       & 900         \\
3        & 100               & 200       & 1200        \\
4        & 150               & 300       & 1800        \\ \hline
District & Demand per period & UAV costs & Truck costs \\ \hline
1        & 200               & 50        & 300         \\
2        & 300               & 100       & 600         \\
3        & 100               & 150       & 900         \\
4        & 150               & 200       & 1200        \\
5        & 250               & 250       & 1500        \\ \hline
District & Demand per period & UAV costs & Truck costs \\ \hline
1        & 200               & 50        & 300         \\
2        & 300               & 100       & 600         \\
3        & 100               & 150       & 900         \\
4        & 150               & 200       & 1200        \\
5        & 250               & 250       & 1500        \\
6        & 200               & 300       & 1800        \\ \hline
\end{tabular}
\end{table}

\begin{table}[ht!]
    \centering
    \footnotesize
    \caption{Demand per district per period for the increasing demand scenario and decreasing demand scenario}
    \label{tab:app:demand_patterns_variants}
\begin{tabular}{llllp{0.1cm}lll}
\hline
       & \multicolumn{3}{l}{Decreasing demand} & & \multicolumn{3}{l}{Increasing demand} \\ \cline{2-4} \cline{6-8} 
Period & District 1  & District 2 & District 3 & & District 1  & District 2 & District 3 \\
\hline
1      & 300         & 450        & 150        & & 100         & 150        & 50         \\
2      & 298         & 447        & 149        & & 102         & 153        & 51         \\
3      & 296         & 444        & 148        & & 104         & 156        & 52         \\
4      & 293         & 440        & 147        & & 107         & 161        & 54         \\
5      & 290         & 435        & 145        & & 110         & 165        & 55         \\
6      & 285         & 428        & 143        & & 115         & 173        & 58         \\
7      & 280         & 420        & 140        & & 120         & 180        & 60         \\
8      & 273         & 410        & 137        & & 127         & 191        & 64         \\
9      & 266         & 399        & 133        & & 134         & 201        & 67         \\
10     & 258         & 387        & 129        & & 142         & 213        & 71         \\
11     & 249         & 374        & 125        & & 151         & 227        & 76         \\
12     & 239         & 359        & 120        & & 161         & 242        & 81         \\
13     & 228         & 342        & 114        & & 172         & 258        & 86         \\
14     & 217         & 326        & 109        & & 183         & 275        & 92         \\
15     & 206         & 309        & 103        & & 194         & 291        & 97         \\
16     & 194         & 291        & 97         & & 206         & 309        & 103        \\
17     & 183         & 275        & 92         & & 217         & 326        & 109        \\
18     & 172         & 258        & 86         & & 228         & 342        & 114        \\
19     & 161         & 242        & 81         & & 239         & 359        & 120        \\
20     & 151         & 227        & 76         & & 249         & 374        & 125        \\
21     & 142         & 213        & 71         & & 258         & 387        & 129        \\
22     & 134         & 201        & 67         & & 266         & 399        & 133        \\
23     & 127         & 191        & 64         & & 273         & 410        & 137        \\
24     & 120         & 180        & 60         & & 280         & 420        & 140        \\
25     & 115         & 173        & 58         & & 285         & 428        & 143        \\
26     & 110         & 165        & 55         & & 290         & 435        & 145        \\
27     & 107         & 161        & 54         & & 293         & 440        & 147        \\
28     & 104         & 156        & 52         & & 296         & 444        & 148        \\
29     & 102         & 153        & 51         & & 298         & 447        & 149        \\
30     & 100         & 150        & 50         & & 300         & 450        & 150        \\
\hline
\end{tabular}
\end{table}

\begin{table}[ht!]
    \centering
    \footnotesize
    \caption{Demand and vehicle allocation costs per district for the Nepal instance}
    \label{tab:app:details_Nepal_instance}
\begin{tabular}{lllll}
\hline
District & Name           & Demand per period & UAV costs & Truck costs \\ \hline
1        & Dolakha        & 217               & 202       & 1256        \\
2        & Gorkha         & 305               & 178       & 1266        \\
3        & Okhaldhunga    & 55                & 266       & 1223        \\
4        & Sindhupalchok  & 278               & 108       & 667         \\
5        & Bhaktapur      & 117               & 26        & 169         \\
6        & Rasuwa         & 49                & 108       & 1928        \\
7        & Ramechhap      & 167               & 214       & 871         \\
8        & Makwanpur      & 156               & 197       & 1085        \\
9        & Dhading        & 352               & 113       & 731         \\
10       & Sindhuli       & 156               & 82        & 683         \\
11       & Nuwakot        & 333               & 67        & 437         \\
12       & Kavrepalanchok & 308               & 62        & 365         \\
13       & Lalitpur       & 107               & 12        & 251         \\
\hline
\end{tabular}
\end{table}

\clearpage

\section{Experimental Settings}
\label{app:settings}
\begin{table}[ht!]
    \centering
    \footnotesize
    \caption{Parameter settings for all methods}
    \label{tab:app:param_settings}
\begin{tabular}{ll}
\hline
General parameters                                            & Setting        \\
\hline
Solving time limit per epoch - Re-optimization        & 900s           \\
Training/solving time limit - other methods           & 14,400s         \\
Shortage margin                                       & $2 \sigma d_{t,t+1,n}$ \\
Algorithmic discount factor $\gamma$                  & 0.9            \\
\hline
Parameters PPO                                        &                \\
\hline
Buffer size                                           & 100            \\
Network learning rate                                 & 0.0001         \\
Network training batch size                           & 64             \\
Network architecture actor and critic                 & 64-64-64       \\
\hline
Parameters DL-VFA                                     &                \\
\hline
Initial $\epsilon$                                    & 0.2            \\
$\epsilon$ decay                                      & 0.98           \\
Initial $\alpha$                                      & 0.2            \\
$\alpha$ decay                                        & 0.99           \\
Buffer size                                           & 1000           \\
Update frequency                                      & 10             \\
\hline
Parameters NN-VFA                                     &                \\
\hline
Initial $\epsilon$                                    & 0.2            \\
$\epsilon$ decay                                      & 0.98           \\
Buffer size                                           & 1000           \\
Update frequency                                      & 10             \\
Network learning rate                                 & 0.001          \\
Network training batch size                           & 256            \\ 
Network architecture                                  & 16-16          \\ \hline
\end{tabular}
\end{table}

\clearpage

\section{Learning Curves}
\label{app:learning_curves}
\begin{figure}[ht!]
    \subfloat[1 district]{\includegraphics[width=0.5\linewidth]{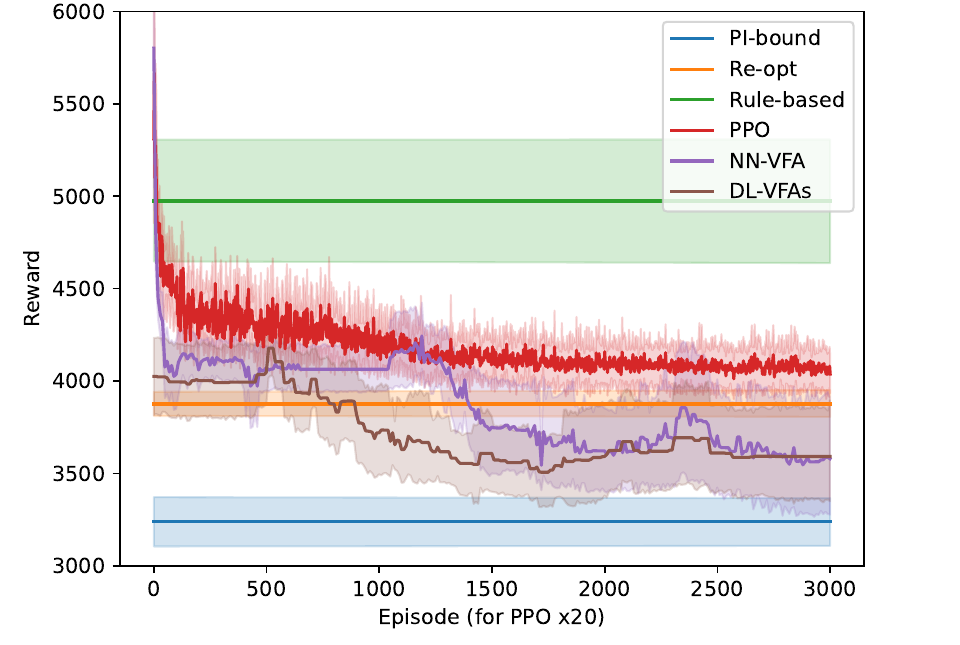}}
    \subfloat[2 districts]{\includegraphics[width=0.5\linewidth]{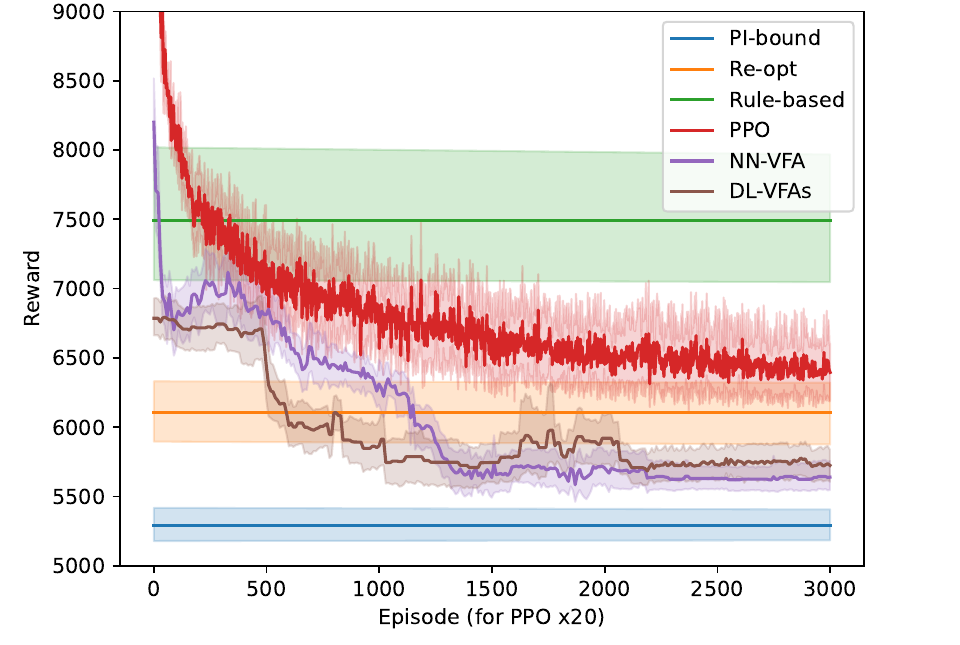}} \\
    \subfloat[4 districts]{\includegraphics[width=0.5\linewidth]{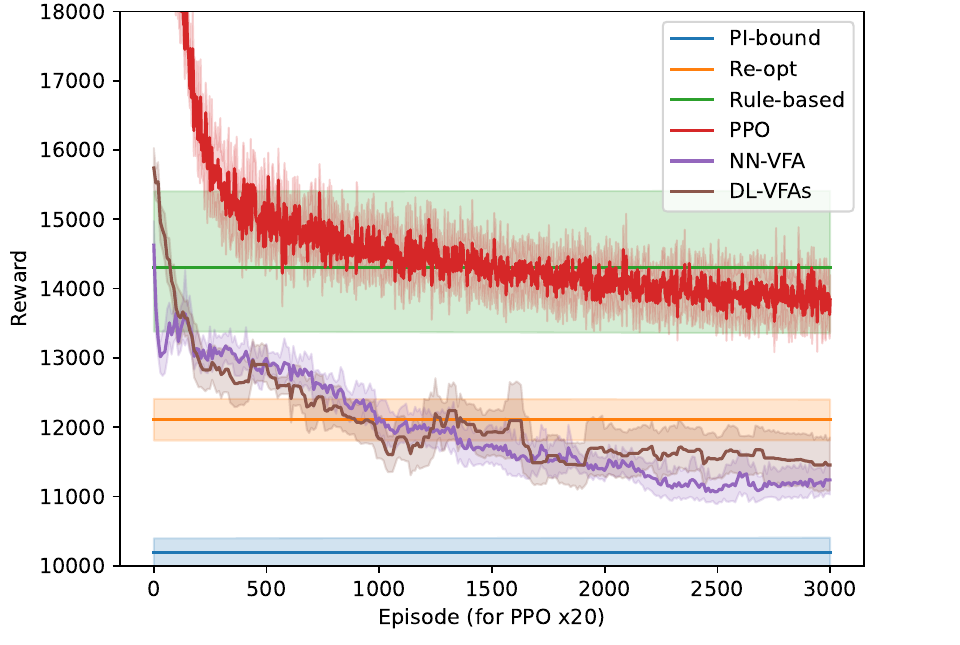}}
    \subfloat[5 districts]{\includegraphics[width=0.5\linewidth]{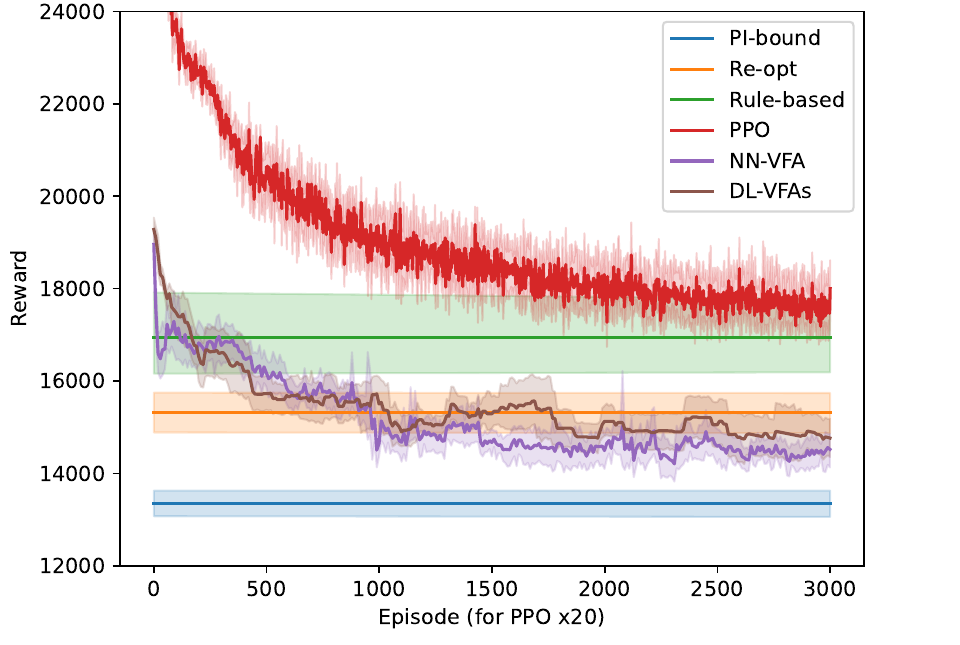}}
    \caption{Learning curves for each algorithm for 1, 2, 4, and 5 districts. The shaded regions around the learning curves provide the 95\% confidence intervals. Each episode is one horizon of 30 periods. For PPO, 60,000 episodes are performed and episodes in the figures are scaled to consist of 20 horizons.}
    \label{fig:app_learning_curves}
\end{figure}

\end{document}